\documentclass[review]{elsarticle}
\RequirePackage[OT1]{fontenc}

\usepackage{amsthm,amsmath,amssymb,mathrsfs,graphicx}
\usepackage{verbatim,multirow,subfigure}
\usepackage{enumerate, enumitem}
\usepackage{color}
\usepackage{array}
\RequirePackage{hyperref}
\usepackage[margin = 1in]{geometry}
\usepackage{setspace}

\usepackage{natbib}
\usepackage{soul}
\theoremstyle{plain}

\newtheorem{theorem}{Theorem}[section]

\newcommand{\msum}{\mathop{\sum}}

\newcommand{\xbar}{\overline{\mathbf{X}}}
\newcommand{\ybar}{\overline{\mathbf{Y}}}

\newcommand{\mtxt}[1]{\mathcal{#1}_n}
\newcommand{\varmn}{{\rm var} \left(\mtxt{M}\right)}

\newcommand{\tr}{{\rm tr}}

\newcommand{\muvec}{\boldsymbol{\mu}}
\newcommand{\pivec}{\boldsymbol{\pi}}
\newcommand{\thetavec}{\boldsymbol{\theta}}
\newcommand{\xvec}{\mathbf{X}}
\newcommand{\yvec}{\mathbf{Y}}
\newcommand{\iid}{{\it i.i.d. \,}}
\newcommand{\etal}{{\it et al. }}
\newcommand{\cov}{{\rm cov}}
\newcommand{\placeref}{{\color{red}place ref }}

\allowdisplaybreaks

\begin{document}
\begin{frontmatter}
	
	\title{High dimensional statistical inference: theoretical development to data analytics}
	
	\author{Deepak Nag Ayyala}
	\address{Department of Population Health Sciences, Medical College of Georgia, Augusta University \\ Augusta, Georgia 30912}

	\begin{abstract}
	This article is due to appear in the Handbook of Statistics, Vol. 43, Elsevier/North-Holland, Amsterdam, edited by Arni S. R. Srinivasa Rao and C. R. Rao.
	
	In modern day analytics, there is ever growing need to develop statistical models to study high dimensional data. Between dimension reduction, asymptotics-driven methods and random projection based methods, there are several approaches developed so far. For high dimensional parametric models, estimation and hypothesis testing for mean and covariance matrices have been extensively studied. However, practical implementation of these methods are fairly limited and are primarily restricted to researchers involved in high dimensional inference. With several applied fields such as genomics, metagenomics and social networking, high dimensional inference is a key component of big data analytics. In this chapter, a comprehensive overview of high dimensional inference and its applications in data analytics is provided. Key theoretical developments and computational tools are presented, giving readers an in-depth understanding of challenges in big data analysis.
	\end{abstract}

	
\end{frontmatter}

\tableofcontents
\section{Introduction}
High dimensional inference and big data analytics are gaining significant prominence in several applied fields such as genomics, imaging neuroscience, econometrics, astronomy and cyber-security \citep{Fan2014}. With accelerated development of technology to study various biological processes and natural phenomenon, there is an exponential growth in the amount of data being generated. Publicly available data sets for genomics such as the Cancer Genome Atlas\footnote{\href{https://portal.gdc.cancer.gov/}{https://portal.gdc.cancer.gov/}} have massive amounts of data that are on the scale of petabytes (1 petabyte = 1024 terabytes). In terms of the number of variables collected (usually represented by $p$) and the number of samples or replicates (represented by $n$), big data can be broadly classified into two categories: (i) large $n$ data sets (ii) large $p$ small $n$ data sets. In data sets with large number of samples, typically arising in astronomy, challenges are mainly computational rather than statistical. Statistical problems in these data sets involve identifying an extremely small number of signals from a large number of observations, {\it a.k.a.} needle in a haystack problem. The large $p$ small $n$ paradigm is commonly encountered in biomedical research areas such as genomics, metagenomics and neuroimaging. The goal in these data sets is to draw inference on a large number of variables simultaneously using a small number of observations. 

Traditional statistical tools are built on the assumption that there is more {\it known} than {\it unkown}, i.e. $n >p$. When $p > n$, asymptotic properties of estimates for parameters such as mean and variance will no longer be valid. For instance consider $p$ parameters $\theta_1, \ldots, \theta_p$ and $\eta = \sum_{k=1}^p \theta_k$ be our parameter of interest. Let $\widehat{\theta}_{nk}$ be a first-order consistent for $\theta_k$ for $k = 1, \ldots, p$, i.e. $\widehat{\theta}_{nk} - \theta_k = o_p(n^{-1/2})$. When $p$ is fixed and finite,  $\widehat{\eta} = \sum_{k=1}^p \theta_{nk}$ will be first-order consistent for $\eta$ because
\begin{equation*}
\widehat{\eta} - \eta = \sum \limits_{k=1}^p (\widehat{\theta}_{nk} - \theta_{k}) = p o_p(n^{-1/2}) = o_p(n^{-1/2}).
\end{equation*}
But if the dimension is a linear function of the sample size, i.e. $p = O(n)$, then this consistency fails because the infinite sum of errors will diverge,
\begin{equation*}
\widehat{\eta} - \eta = p o(n^{-1/2}) = o_p(n^{1/2})
\end{equation*}
This problem can be solved by considering second-order consistent estimators, which require additional attention to the asymptotic properties of $\widehat{\theta}_{nk}$ to derive. In multivariate models, there are two main parameters of interest - mean vector and variance matrix. Estimation of these parameters and construction of hypothesis tests for high dimensional data require additional calculations to have good asymptotic behaviour. 

Large data sets with discrete data are very commonly observed in various fields. In text mining, the distribution of words in a document are recorded by counting the number of occurrences of each word in the document. In genomics and metagenomics, recently developed high-throughput experimental procedures are making it possible to record counts of genes expression and bacterial abundance in samples. However, statistical literature on multivariate models for discrete data is very sparse. Most of the multivariate probability models that are encountered in high dimensional literature are continuous. Unlike the continuous distributions, multivariate analogues of standard discrete models such as Bernoulli, binomial, Poisson, etc. are not extensively developed. In the univariate case, mixture models such as beta-binomial and Poisson-gamma have been developed to address over-dispersion in count data. Multivariate models do not exist for all mixture distributions. 

In this chapter, we will look at three topics of interest in high dimensional inference. In Section \ref{sec:meanvt}, we will look at hypothesis tests for the mean vector. Estimation and hypothesis tests for the covariance matrix will be addressed in Section \ref{sec:covariancematrix}. Formulation of standard discrete multivariate models and parameter estimation of hierarchical multivariate count models are presented in Section \ref{sec:discretemodels}. Finally, we will conclude with some challenges that still lie ahead of us in high dimensional inference in Section \ref{sec:conclusion}.

\section{Mean vector testing}
\label{sec:meanvt}
The first moment, mean, is the most commonly studied parameter when exploring the properties of distributions. The mean or expected value of a random variable is a measure of location of the center of the data. When comparing $p$ dimensional two distributions, equality of means indicates that distributions are centered around the same point in the sample space. Given two variables characterized by their means, $\mathbf{X} \sim \mathcal{F}_p (\cdot, \muvec_1)$, and $\mathbf{Y} \sim \mathcal{G}_p (\cdot, \muvec_2)$, the hypothesis of comparing means can be stated as
\begin{equation}
H_0 : \muvec_1 = \muvec_2 \hspace{5mm} \mbox{vs. } \hspace{5mm} H_A : \muvec_1 \neq \muvec_2.
\label{eqn:meanhyp1}
\end{equation}
Given samples $\xvec_1, \ldots, \xvec_n$ and $\yvec_1, \ldots, \yvec_m$ from the two distributions, sample means $\xbar = n^{-1} \msum_{i = 1}^n \xvec_i$ and $\ybar = \msum_{j = 1}^m \yvec_j$ are natural unbiased estimators of $\muvec_1$ and $\muvec_2$ respectively. Hence the difference of sample averages $\overline{X} - \overline{Y}$ will be unbiased for $\muvec_1 - \muvec_2$. To calculate a test statistic, a functional needs to be defined to map the multivariate difference $\xbar - \ybar$ on to the real line. 

Hotelling's $T^2$ \citep{Hotelling1931} was the first such test constructed which uses the Mahalanobis distance as the functional,
\begin{equation}
T^2_{Hot} = \frac{ n + m - p - 1}{(n + m - 2)p} \, \frac{n m}{n + m} \left( \xbar - \ybar \right)^{\top} \mathcal{S}^{-1} \left( \xbar - \ybar \right)
\label{eqn:hotelingt2}
\end{equation} 
where $\mathcal{S} = (n + m - 2)^{-1} \left\{ \msum_{i = 1}^n \left( \xvec_i - \xbar \right) \left( \xvec_i - \xbar \right)^{\top} + \msum_{j = 1}^m \left( \yvec_j - \ybar \right) \left( \yvec_j - \ybar \right)^{\top}  \right\}$ is the pooled sample covariance matrix with ${\rm rank}(\mathcal{S}) = \min(p, n + m - 2)$. Under $H_0$, the test statistic follows a $F_{p, n + m - p - 1}$ distribution provided $\mathcal{F}_p$ and $\mathcal{G}_p$ are both homogeneous multivariate Gaussian distributions with common covariance matrix $\Sigma$ and $p < n + m - 1$. 

The Hotelling's $T^2$ test is developed for the two-sided alternative in \eqref{eqn:meanhyp1}. The functional in $T^2_{Hot}$ has a quadratic form and is always non-negative as $\mathcal{S}$ is positive definite. Hence $T^2_{Hot}$ does not differentiate between $\muvec_1 - \muvec_2 = \boldsymbol{\delta}$ and $\muvec_1 - \muvec_2 = - \boldsymbol{\delta}$. To define a one-sided alternative, an order in $\mathbb{R}^p$ should first be determined. For example, the lexicographic order or the partial element-wise order can be considered. For the one sample case, Kudo \citep{Kudo1963} developed a likelihood ratio test (LRT) for the alternative $H_A : \muvec_1 > \mathbf{0}$ where the inequality indicates $\mu_{1i} \geq 0$ for $i = 1, \ldots, p$ with at least one $\mu_{1i} > 0$. Maximizing the parameters over the positive cone is done using quadratic programming. p-value calculation is computationally intensive due to the $2^p$ potential pairs of zeros and positive elements under the alternative. Also a two-sample extension for this test was not provided.

The Hotelling's $T^2$ test has three major deficiencies when doing inference in high dimensions. 
\begin{enumerate}[label=\underline{Issue \Roman*}]
	\item The test is defined only when $p < n + m - 2$. Typically in data sets arising in genomics and other high throughput experiments, the dimension is in thousands and the number of samples are in tens or hundreds.  
	\item The test holds only for comparing Gaussian distributions, an assumption that is not straightforward to verify in higher dimensions. In distributions with restricted sample space such as the Dirichlet, the mean vector is a location parameter.  
	\item The test requires observations to be {\it independently and identically distributed}, i.e. \iid In imaging studies such as fMRI experiments, the observations are not \iid The inherent dependence structure in such data sets is ignored by $T^2_{Hot}$, potentially leading to biased estimates.
\end{enumerate}
To address these shortcomings, one has to use test statistics that take into account the bias due to high dimension and the dependence structure in the data. For instance to address Issue I, there are two approaches that can used. The first method is to study the asymptotic properties of a functional of $\xbar - \ybar$ and construct a large-sample test. This method can also relax Issue II by accommodating non-Gaussian distributions through conditions on the moments of the distribution. Second method is to reduce dimension by projecting the $p$-variate samples into a lower dimensional space such that traditional tests such as $T^2_{Hot}$ can be applied. The dependence structure in Issue III is complicated since the entire autocovariance function of the model needs to be considered. Parametrizing the autocovariance function and restricting the dependence structure can help reduce the complexity of the problem. 

\subsection{Independent observations}
\label{sec:asymptoticmvt}
First let us address testing the hypothesis in \eqref{eqn:meanhyp1} for \iid samples in high dimension. When $p > n + m - 2$, the pooled sample covariance matrix $\mathcal{S}$ is rank-deficient and does not have a well-defined inverse. The Mahalanobis distance of $\muvec_1 - \muvec_2$ is not a valid measure. To construct a test statistic, we need a functional of $\xbar - \ybar$ which is zero in expectation when $H_0$ is true and non-zero when $H_A$ is true. A Natural choice of such functional which does not involve $\mathcal{S}$ is the $\ell_d$-norm for $d > 0$. When $d = 1$, Chung and Fraser \citep{Chung1958} proposed a permutation test using the sum of element-wise $t$-test statistics,
\begin{equation}
T_{CF} = \msum_{k = 1}^p \frac{ \left| \overline{X}_k - \overline{Y}_k \right|}{\mathcal{S}_{kk}},
\label{eqn:chungfraser}
\end{equation}
as the test statistic. The Euclidean norm\footnote{Abuse of notation: The squared Euclidean norm is referred to as the Euclidean norm unless otherwise stated} is preferred over the $\ell_1$-norm due to ease of calculation of moments. Dempster \citep{Dempster1958} developed the first test statistic using the Euclidean norm of difference of means, $\left( \xbar - \ybar \right)^{\top} \left( \xbar - \ybar \right)$. The test statistic is given by
\begin{equation}
T_{Demp} = \frac{ \left( \xbar - \ybar \right)^{\top} \left( \xbar - \ybar \right)}{ \msum_{k = 1}^{n + m - 2} \mathbf{W}_k^{\top} \mathbf{W}_k },
\label{eqn:dempster}
\end{equation}
where $\{\mathbf{W}_k, k = 1, \ldots, n + m - 2 \}$ are orthogonal vectors such that the set of vectors $\{(n + m)^{-1} (n \xbar + m \ybar), \xbar - \ybar, \mathbf{W}_1, \ldots, \mathbf{W}_{n + m -2} \}$ form an orthogonal basis for the space spanned by $\{\xvec_1, \ldots, \xvec_n, \yvec_1, \ldots, \yvec_m \}$. The Dempster test is {\it non-exact} and is distributed as an $F_{r, (n + m - 2)r}$ under the null hypothesis. The parameter $r$ is unknown and is estimated from the data. However, both these tests ignore the covariance structure and are shown to not perform well even when $p$ is close to $n + m$.

To construct a large sample test, the asymptotic properties of the Euclidean norm of $\xbar - \ybar$ need to be studied. When the two distributions are homogeneous with covariance matrix $\Sigma$, we have
\begin{equation}
\mathbb{E} \left\{  \left( \xbar - \ybar \right)^{\top} \left( \xbar - \ybar \right) \right\} = \left( \muvec_1 - \muvec_2 \right)^{\top} \left( \muvec_1 - \muvec_2 \right) + \left(\frac{1}{n} + \frac{1}{m} \right) \tr \left( \Sigma \right).
\label{eqn:euclideanmean}
\end{equation}
Without loss of generality, assume $n < m$. Under $H_0$,  $\left( \xbar - \ybar \right)^{\top} \left( \xbar - \ybar \right)$ has expected value equal to $\mathcal{B}_n = 2 (1/n + 1/m) \tr \left( \Sigma \right)  \leq 2 n^{-1} p \, \lambda_{\max}$, where $\lambda_{\max}$ is the largest eigenvalue of $\Sigma$. If $p$ is fixed, then $ \lim\limits_{ n \rightarrow \infty} \mathcal{B}_n = 0$, implying $ \left( \xbar - \ybar \right)^{\top} \left( \xbar - \ybar \right)$ is asymptotically unbiased. But if $p$ increases with $n$, then the Euclidean norm needs to be adjusted for this bias. For instance, if we assume $p = Cn^{\alpha}$ for some $\alpha > 0$, then $\mathcal{B}_ n = 2 n^{1 - \alpha} \lambda_{\max}$ which diverges when $\alpha > 1$. Further note that properties of $\mathcal{B}_n$ are independent of the distributions of the two groups.

To adjust the bias, consider the pooled sample covariance matrix $\mathcal{S}$, which is unbiased for $\Sigma$. Since trace is a linear functional, $\tr (\mathcal{S})$ will be unbiased for $\tr (\Sigma)$. This gives
\begin{equation*}
\mathcal{M}_n = \left( \xbar - \ybar \right)^{\top} \left( \xbar - \ybar \right) - \frac{n + m}{nm} \tr \left( \mathcal{S} \right).
\end{equation*}
 as an unbiased estimator of $ \left( \muvec_1 - \muvec_2 \right)^{\top} \left( \muvec_1 - \muvec_2 \right)$. Using its quadratic form, the variance of $\mathcal{M}_n$ can be calculated as $\varmn = 2 (1/n + 1/m)^2 \{1 + 1/(n + m - 2) \} \tr \Sigma^2 \{1 + o(1) \}$. The error term, $1 + o(1)$, vanishes under Gaussian assumption. To construct a test statistic using $\mathcal{M}_n$, a ratio consistent estimator of $\varmn$ is needed. 
 
  In their seminal work, Bai and Saranadasa \citep{Bai1996} used $\mathcal{M}_n$ to construct the test statistic
 \begin{equation}
 T_{BS} = \frac{ \left( \xbar - \ybar \right)^{\top} \left( \xbar - \ybar \right) - \frac{n + m}{nm} \tr \left( \mathcal{S} \right)}{ \frac{n + m}{nm} \sqrt{ \frac{ 2(n + m - 1)(n + m - 2)}{(n + m)(n + m -3)} \left\{ \tr \left(\mathcal{S}^2\right) - (n + m - 2)^{-1} \tr^2 \mathcal{S}  \right\} }}
 \label{eqn:baisaranadasa}
 \end{equation}
The test statistic follows a standard normal distribution asymptotically under the following conditions:
\begin{enumerate}[label=(BS \Roman*)]
	\item \label{itm:bs1} $p/n  \rightarrow \delta > 0$, indicating that $p$ increases faster than $n$.
	\item \label{itm:bs2} $n/(n + m) \rightarrow \kappa \in (0,1)$ meaning sample sizes from both groups have proportionate rates of increase.
	\item \label{itm:bs3} $\lambda_{\max} = o \left( \sqrt{ \tr \Sigma^2} \right)$, which relates to the strength of the covariance structure. 
	\item \label{itm:bs4} $(\muvec_1 - \muvec_2)^{\top} \Sigma (\muvec_1 - \muvec_2) = o \left( (1/n + 1/m) \, \tr \Sigma^2 \right)$ is a local alternative condition to calculate the asymptotic power, under which the variance estimate remains ratio consistent.
\end{enumerate}

Let us elaborate condition \ref{itm:bs3} to understand how strong the covariance structure can be. Consider the independent elements case, $\Sigma = \mathcal{I}$ with $\lambda_{\max} = 1$ and $\tr \Sigma^2 = p$. Thus we have $\lambda/\sqrt{ \tr \Sigma^2} = 1/\sqrt{p} \rightarrow 0$, which indicates the validity of the condition. If we consider a moving average covariance structure with $\Sigma_{ij} = \rho^{|i - j|}$ for $0 < \rho < 1$. Then we have $\lambda_{\max} \leq (1 + \rho)/(1 - \rho)$ and $\tr \Sigma^2 \approx p (1 - \rho^p)(1 - \rho)^{-1}$ which also satisfies the condition for all values of $\rho$. The condition, however, does not allow covariance structures from the other end of the spectrum: an exchangeable covariance structure with $\Sigma_{ij} = \rho$ for all $i,j$ and for some $0 < \rho < 1$ which has $\lambda_{\max} = 1 + (p -  1)\rho$ and $\tr \Sigma^2 = p + (p^2 - p) \rho^2$. This gives 
\begin{equation*}
\lim \limits_{p \rightarrow \infty} \frac{ \lambda_{\max}}{\sqrt{\tr \Sigma^2}} = \lim \limits_{p \rightarrow \infty} \frac{1 + (p - 1) \rho}{\sqrt{p + (p^2 - p) \rho^2}} = 1,
\end{equation*}
which does not satisfy the condition. 

The Bai-Saranadasa test statistic is highly regarded in high dimensional mean-vector testing literature. In addition to extending the test to higher dimensions, it also relaxed the normality assumption on the samples. Instead, the observations are assumed to be coming from a factor model of the form
\begin{equation}
\xvec = \muvec + \Gamma \mathbf{Z},
\label{eqn:factormodel}
\end{equation}
where $\mathbf{Z} = (Z_1, \ldots, Z_p)$ and $Z_i$'s are continuous \iid random variables with $E(Z_i) = 0$ and $E(Z_i^4) = 3 + \Delta < \infty$. The covariance structure is determined by $\Gamma$ through the relationship $\Sigma = \Gamma \Gamma^{\top}$. When $\Delta = 0$, the elements of $\mathbf{Z}$ are normally distributed. When $0 < \Delta < \infty$, the $Z_i$'s have heavier tails than normal, yet have finite moments. Examples of distributions satisfying the moment conditions are Laplace or double exponential distribution and centered gamma distribution. 

In equation \eqref{eqn:euclideanmean}, the trace term comes only from the inner products of $\xvec_i$'s and $\yvec_j$'s. For any $i$, we have $\mathbb{E}(\xvec_i^{\top} \xvec_i) = \muvec_1^{\top} \muvec_1+ \tr n^{-1} \Sigma$ and $\mathbf{E}(\xvec_i^{\top} \xvec_j^{\top}) = \muvec_1^{\top} \muvec_1$ when $i \neq j$. Hence subtracting the inner-product terms from $n^2 \mathbf{E}(\overline{\xvec}^{\top} \overline{\xvec})$ and $m^2 \mathbf{E}(\overline{\yvec}^{\top} \overline{\yvec})$, we have
\begin{equation*}
\mathbf{E} \left( \msum_{i \neq j}^n \xvec_i^{\top} \xvec_j \right) = n (n - 1) \muvec_1^{\top} \muvec_1, \,\, 
\mathbf{E} \left( \msum_{i \neq j}^m \yvec_i^{\top} \yvec_j \right) = m (m - 1) \muvec_2^{\top} \muvec_2, \,\, 
\mathbf{E} \left( \msum_{i, j} \xvec_i^{\top} \yvec_j \right) = n m \muvec_1^{\top} \muvec_2.
\end{equation*}
Combining the terms in the above equation, the statistic
\begin{equation}
\mathcal{T}_n = \frac{1}{n(n-1)} \msum_{i \neq j}^n \xvec_i^{\top} \xvec_j + \frac{1}{m(m-1)} \msum_{i \neq j}^m \yvec_i^{\top} \yvec_j - \frac{2}{nm} \msum_{i = 1}^n \msum_{j = 1}^m \xvec_i^{\top} \yvec_j 
\label{eqn:euclideanmean2}
\end{equation}
has expected value equal to $(\muvec_1 - \muvec_2)^{\top} (\muvec_1 - \muvec_2)$. 

Chen and Qin \citep{Chen2010} constructed a test statistic using $\mathcal{T}_n$ as the functional, which has zero expected value under $H_0$. They assumed that the data follows the factor model in equation \eqref{eqn:factormodel}. Sample sizes are restricted similar to \ref{itm:bs2}. A major criticism of $T_{BS}$ has been the restriction of homogeneity of the two populations, i.e. equal covariance structure. Addressing this issue is a major achievement of the Chen and Qin test, which relaxed this condition. The two populations are allowed to have unequal covariance structures, $\Sigma_1$ and $\Sigma_2$ respectively. This extension results in the local alternative condition in \ref{itm:bs4} to be modified, with the rate holding with respect to both $\Sigma_1$ and $\Sigma_2$. Strength of the covariance matrix as restricted by \ref{itm:bs3} is also modified to accommodate the heterogeneity. Another major accomplishment of the Chen and Qin test is removing a direct constraint between $p$ and $n$ as in \ref{itm:bs1}.

The modified constraints on the model are summarized as follows:
\begin{enumerate}[label=(CQ \Roman*)]
\addtocounter{enumi}{2}
	\item \label{itm:cq3} $\tr \left( \Sigma_a \Sigma_b \Sigma_c \Sigma_d \right) = o \left[ \tr^2 \left\{ \left(\Sigma_1 + \Sigma_2 \right)^2  \right\}  \right]$ for $a, b, c, d \in \{1, 2\}$.  
	\item \label{itm:cq4} $\left(\muvec_1 - \muvec_2 \right)^{\top} \Sigma_a \left( \muvec_1 - \muvec_2 \right) = o \left[ (n + m - 2)^{-1} \tr \left\{ \left( \Sigma_1 + \Sigma_2 \right)^2 \right\} \right]$ for $a = 1, 2$.
\end{enumerate}
Under the local alternative, variance of $\mathcal{T}_n$ is equal to 
\begin{equation*}
{\rm var}\left(\mathcal{T}_n \right) = \left[  \frac{2}{n(n-1)} \tr\left( \Sigma_1^2 \right) + \frac{2}{m(m-1)} \tr \left( \Sigma_2^2 \right) + \frac{4}{nm} \tr \left( \Sigma_1 \Sigma_2 \right) \right] \left\{ 1 + o(1) \right\}.
\end{equation*}

As used in $T_{BS}$, $\left\{ \tr (\mathcal{S}_1^2) - n^{-1} \tr^2 \mathcal{S}_1 \right\}$ can be used as a ratio consistent estimator of $\tr (\Sigma_1^2)$. Inspired by the removal of inner-product terms in $\mathcal{T}_n$, Chen and Qin argue that similar rationale relaxes a direct relationship between $p$ and $n$ as in \ref{itm:bs1}.  They proposed ratio consistent estimators of the form
\begin{align*}
\widehat{\tr (\Sigma_1^2)} &= \frac{1}{n(n-1)} \tr \left\{ \msum_{i = 1}^n \msum_{j \neq i}\left( \xvec_i - \overline{\xvec}_{(i,j)}  \right) \xvec_i^{\top} \left( \xvec_j - \overline{\xvec}_{(i,j)} \right) \xvec_j^{\top} \right\}, \\
\widehat{\tr (\Sigma_2^2)} &= \frac{1}{m(m-1)}  \tr \left\{ \msum_{i = 1}^m \msum_{j \neq i} \left( \yvec_i - \overline{\yvec}_{(i,j)} \right) \yvec_i^{\top} \left( \yvec_j - \overline{\yvec}_{(i,j)} \right) \yvec_j^{\top} \right\}, \\
\widehat{\tr (\Sigma_1 \Sigma_2)} &= \frac{1}{n m} \tr \left\{  \msum_{i = 1}^n \msum_{j =1}^m \left( \xvec_i - \overline{\xvec}_{(i)} \right) \xvec_i^{\top} \left( \yvec_j - \overline{\yvec}_{(j)} \right) \yvec_j^{\top} \right\},
\end{align*}
where $\overline{\xvec}_{(i)} = (n - 1)^{-1} \msum_{k \neq i}^n \xvec_k$, $\overline{\xvec}_{(i,j)} = (n - 2)^{-1} \msum_{k \neq i, j}^n \xvec_k$, $\overline{\yvec}_{(i)} = (n - 1)^{-1} \msum_{k \neq i}^n \yvec_k$ and $\overline{\yvec}_{(i,j)} = (n - 2)^{-1} \msum_{k \neq i, j}^n \yvec_k$. Finally, the Chen-Qin test statistic is given by
\begin{equation}
T_{CQ} = \frac{ \mathcal{T}_n}{ \sqrt{\frac{2}{n(n-1)} \widehat{ \tr\left( \Sigma_1^2 \right)} + \frac{2}{m(m-1)} \widehat{\tr \left( \Sigma_2^2 \right)} + \frac{4}{nm} \widehat{\tr \left( \Sigma_1 \Sigma_2 \right)} }}
\label{eqn:chenqin}
\end{equation}
which follows a normal distribution asymptotically under $H_0$. 

In $T_{BS}$ and $T_{CQ}$, the Euclidean norm is used as the functional to avoid inverting the sample covariance matrix, which is singular when $p > n$. While $\mathcal{S}$ is not invertible, the diagonal elements are all non-zeroes and invertible (a zero diagonal element implies the corresponding variable is a constant and it can be removed from the analysis). Using the diagonal elements, a {\it modified} Mahalanobis distance can be calculated as a weighted Euclidean norm,
\begin{equation*}
\mathcal{W}_n = \left( \overline{\xvec} - \overline{\yvec} \right)^{\top} \mathcal{D}_{\mathcal{S}}^{-1} \left( \overline{\xvec} - \overline{\yvec} \right) = \msum_{k = 1}^p \frac{ \left( \overline{X}_k - \overline{Y}_k \right)^2}{\mathcal{S}_{kk}},
\end{equation*}
where $\mathcal{D}_{\mathcal{S}}$ is the $p \times p$ diagonal matrix of $\mathcal{S}$. When the two groups are homogeneous, we have $\mathbb{E} \left( \overline{X}_k - \overline{Y}_k \right)^2  = \left( \mu_{1k} - \mu_{2k} \right)^2 + \left(1/n + 1/m \right) \sigma_{kk}$ and $\mathbb{E} \left( \mathcal{S}_{kk} \right) = (n + m - 2)/(n + m) \, \sigma_{kk}$. As the ratio of these two expected values independent of the index $k$, we have
\begin{equation*}
\mathbb{E} \left( \mathcal{W}_n \right) =  \left( \muvec_1 - \muvec_2 \right)^{\top} \left( \muvec_1 - \muvec_2 \right) + \left( \frac{1}{n} + \frac{1}{m} \right) \frac{ n + m}{n + m - 2}\, p.
\end{equation*}
Similar to the calculations in the Euclidean norm, it is straightforward to show using the quadratic form that ${\rm var} \left( \mathcal{W}_n \right) = 2 \tr \left( \mathcal{R}^2 \right) \left\{1 + o(1) \right\}$.

Srivastava and Du \citep{SrivastavaDu2008} developed a test statistic based on $\mathcal{W}_n$ as the functional, adjusting for its expected value. The test statistic is valid under the following assumptions:
\begin{enumerate}[label=(SD \Roman*)]
	\item \label{itm:sd1} The dimension increases at a polynomial rate with respect to $n$, $n = O \left( p^{\delta} \right), 1/2 < \delta \leq 1$.
	\item \label{itm:sd2} Sample sizes of the two groups, $n$ and $m$, are constrained as in \ref{itm:bs2}.
	\item \label{itm:sd3} If $\mathcal{R}$ is the population correlation matrix and $\lambda_1 \geq \ldots \geq \lambda_p$ are its eigenvalues, then $\lim_{p \rightarrow \infty} \tr \left(\mathcal{R}^k \right)/p \in \left(0, \infty \right)$ for $k = 1,2,3,4$ and $\lambda_1 = o \left( \sqrt{p} \right)$.
	\item \label{itm:sd4} Means of the two groups satisfy the local alternative condition: $ \left( \muvec_1 - \muvec_2 \right)^{\top} \mathcal{D}_{\Sigma}^{-1} \left(\muvec_1 - \muvec_2 \right) \leq M p/(n + m -2) (1/n + 1/m)$ for some finite constant $M$.
\end{enumerate}
The Srivastava-Du test statistic is given by
\begin{equation}
T_{SD} = \frac{ \frac{n m}{n + m} \left(\overline{\xvec} - \overline{\yvec} \right)^{\top} \mathcal{D}_{\mathcal{S}}^{-1} \left(\overline{\xvec} - \overline{\yvec} \right) - \frac{ (n + m) \, p}{n + m - 2}}{ \sqrt{2 \, \left(\tr R^2 - p^2/n \right) \left( 1 + \tr R^2/p^{3/2} \right)}}
\label{eqn:srivastavadu}
\end{equation}
where $R = \mathcal{D}_{\mathcal{S}}^{-1/2} \mathcal{S} \mathcal{D}_{\mathcal{S}}^{-1/2}$ is the sample correlation matrix. The test statistic is asymptotically normal under the null hypothesis. 

The condition imposed on the correlation structure in \ref{itm:sd3} is very restrictive compared to \ref{itm:bs3} and \ref{itm:cq3}. For example, consider $\Sigma = \mathcal{R} = {\rm diag}(p^{\omega}, 1, \ldots, 1)$ for some $1/4 \leq \omega < 1$. Then $\tr \Sigma^2 = p + p^{2\omega} - 1$, $\tr \Sigma^4 =  p + p^{4 \omega} - 1$ and $\lambda_{\max} = p^{\omega}$. \ref{itm:bs3} and \ref{itm:cq3} are satisfied as 
\begin{equation*}
\frac{\lambda_{\max}}{\tr \Sigma^2} = \frac{ p^{\omega}}{p + p^{2 \omega} - 1} \rightarrow 0, \quad
\frac{ \tr \Sigma^4}{ \tr^2 \Sigma^2} = \frac{ p + p^{4 \omega} - 1}{ \left( p + p^{2 \omega} - 1 \right)^2} \rightarrow 0.
\end{equation*}
For \ref{itm:sd3}, we have $\lambda_{\max}/\sqrt{p} = p^{\omega - 1/2} \rightarrow 0$ but $ \tr \mathcal{R}^4/p = \left( p + p^{4 \omega} - 1 \right)/p = 1 + p^{4 \omega - 1} - p^{-1}$, which is not bounded for $\omega > 1/4$.

Another major constraint of the Srivastava-Du test is that the observations are assumed to be normally distributed. Unlike $T_{BS}$ and $T_{CQ}$, asymptotically equivalent expressions for ${\rm var}(\mathcal{W}_n)$ are not established. Instead, exact variance is derived using the properties of the normal distribution. In a sequence of papers, the authors have provided extensions to $T_{SD}$ to reduce some of the assumptions. In Srivastava \citep{Srivastava2009}, the term $\tr R^2/p^{3/2}$ in the denominator of $T_{SD}$ was shown to converge to zero and hence dropped. In Srivastava-Kano \citep{Srivastava2013}, an extension to the heterogeneous case was developed. However this test is {\it inexact} in the sense that the functional $\mathcal{W}_n^*$ has expected value equal to $\left( \muvec_1 - \muvec_2 \right)^{\top} \left( \muvec_1 - \muvec_2 \right)$ only in limit. 

Inspired by the idea of Chen and Qin \citep{Chen2010}, Park and Ayyala \citep{Park2013} modified the functional $\mathcal{W}_n$ by removing the inner product terms.  Using the true covariance diagonal, the functional
\begin{equation}
\mathcal{U}_n^* =   \frac{1}{n(n - 1)} \msum_{i =\neq j}^n \xvec_i^{\top} \mathcal{D}_{\Sigma}^{-1} \xvec_j + \frac{1}{m(m - 1)} \msum_{i \neq j}^m \yvec_i^{\top} \mathcal{D}_{\Sigma}^{-1} \yvec_j - \frac{2}{n \, m} \msum_{i = 1}^n \msum_{j = 1}^m \xvec_i^{\top} \mathcal{D}_{\Sigma}^{-1} \yvec_j
\end{equation}
has expected value $\left( \muvec_1 - \muvec_2 \right)^{\top} \left( \muvec_1 - \muvec_2 \right)$. Replacing the true covariances with consistent estimators, a leave-out approach has been implemented to maintain independence amongst the terms. For instance in $\xvec_i^{\top} \mathcal{D}_{\Sigma}^{-1} \xvec_j$, the quantities $\xvec_i, \xvec_j$ and $\widehat{\mathcal{D}_{\Sigma}}$ will be independent if $\widehat{\mathcal{D}_{\Sigma}}$ is constructed by {\it leaving out} $\xvec_i$ and $\xvec_j$. The pooled sample covariance matrix $\mathcal{S} = \left( (n - 1) \mathcal{S}_1 + (m - 1) \mathcal{S}_2 \right)/(n + m - 2)$, where $\mathcal{S}_1$ and $\mathcal{S}_2$ are the sample covariance matrices of the two groups respectively, is not useful because $\mathcal{S}_1$ contains $\xvec_i$ and $\xvec_j$. If these two samples are removed from $\mathcal{S}_1$, then $\mathcal{S}_{1(i,j)} = (n - 3)^{-1} \msum_{k \neq i, j}^n \left( \xvec_k - \overline{\xvec}_{(i,j)} \right) \left( \xvec_k - \overline{\xvec}_{(i,j)} \right)^{\top}$, where $\overline{\xvec}_{(i,j)} = (n - 2)^{-1} \msum_{k \neq i,j}^n \xvec_k$ will be independent of $\xvec_i$ and $\xvec_j$.  Similarly, for the second and third terms, we can define $\mathcal{S}_{2(i,j)}$, $\mathcal{S}_{1(i)}$ and $\mathcal{S}_{2(j)}$ respectively to maintain independence of terms. Then diagonals of the pooled sample estimators 
\begin{equation*}
\mathcal{S}^{(1)}_{(i,j)} = \frac{ (n - 3) \mathcal{S}_{1(i,j)} + (m - 1) \mathcal{S}_{2}}{n + m - 4}, \,
\mathcal{S}^{(2)}_{(i,j)} = \frac{ (n - 1) \mathcal{S}_{1} + (m - 3) \mathcal{S}_{2(i,j)}}{n + m - 4}, \,
\mathcal{S}^{(12)}_{(i,j)} = \frac{ (n - 2) \mathcal{S}_{1(i)} + (m - 2) \mathcal{S}_{2(j)}}{n + m - 4},
\end{equation*}
is used  to construct the functional
\begin{equation}
\mathcal{U}_n =  \frac{n + m - 6}{n + m - 4} \left( \frac{1}{n(n - 1)} \msum_{i \neq j}^n \xvec_i^{\top} \mathcal{D}_{\mathcal{S}^{(1)}_{(i,j)}}^{-1} \xvec_j + \frac{1}{m(m - 1)} \msum_{i \neq j}^m \yvec_i^{\top} \mathcal{D}_{\mathcal{S}^{(2)}{(i,j)}}^{-1} \yvec_j - \frac{2}{n \, m} \msum_{i = 1}^n \msum_{j = 1}^m \xvec_i^{\top} \mathcal{D}_{\mathcal{S}^{(12)}_{(i,j)}}^{-1} \yvec_j \right),
\label{eqn:parkayyala1}
\end{equation}
which has expected value $\left( \muvec_1 - \muvec_2 \right)^{\top} \left( \muvec_1 - \muvec_2 \right)$. 

From the quadratic form and independence of the terms, variance of $\mathcal{U}_n$ will be 
\begin{equation*}
{\rm var}\left( \mathcal{U}_n \right) = \left( \frac{n + m - 6}{n + m - 4} \right)^2 \left\{ \frac{2}{n(n - 1)} \tr \left( \mathcal{R}_1^2 \right) + \frac{2}{m(m-1)} \tr \left( \mathcal{R}_2^2  \right) + \frac{4}{nm} \tr \left( \mathcal{R}_1 \mathcal{R}_2  \right)    \right\},
\end{equation*}
where $\mathcal{R}_1$ and $\mathcal{R}_2$ are the for notational  convenience to identify that the terms are related to $\xvec$ and $\yvec$ respectively. A similar {\it leave-out} approach is applied to modify the standard correlation matrix estimate $\widehat{\mathcal{R}_1} = \mathcal{D}_{\mathcal{S}_1}^{-1/2} \mathcal{S}_1 \mathcal{D}_{\mathcal{S}}^{-1/2}$. Centering the observations only once as in $T_{CQ}$ and rearranging the terms, the estimators
\begin{align*}
	\widehat{\tr \left( \mathcal{R}_1^2 \right)} = \frac{1}{n(n-1)} \tr \left\{ \msum_{i = 1}^n \msum_{j \neq i} \xvec_i^{\top} \mathcal{D}_{\mathcal{S}^{(1)}_{(i,j)}}^{-1}  \left( \xvec_j - \overline{\xvec}_{(i,j)} \right)    \xvec_j^{\top} \mathcal{D}_{\mathcal{S}^{(1)}_{(i,j)}}^{-1}  \left( \xvec_i - \overline{\xvec}_{(i,j)} \right)    \right\}, \\
	\widehat{\tr \left( \mathcal{R}_2^2 \right)} = \frac{1}{m(m-1)} \tr \left\{ \msum_{i = 1}^m \msum_{j \neq i} \yvec_i^{\top} \mathcal{D}_{\mathcal{S}^{(2)}_{(i,j)}}^{-1}  \left( \yvec_j - \overline{\yvec}_{(i,j)} \right)    \yvec_j^{\top} \mathcal{D}_{\mathcal{S}^{(2)}_{(i,j)}}^{-1}  \left( \yvec_i - \overline{\yvec}_{(i,j)} \right)    \right\}, \\
	\widehat{\tr \left( \mathcal{R}_{1} \mathcal{R}_2 \right)} = \frac{1}{n m} \tr \left\{ \msum_{i = 1}^n \msum_{j = 1}^m \xvec_i^{\top} \mathcal{D}_{\mathcal{S}^{(12)}_{(i,j)}}^{-1}  \left( \yvec_j - \overline{\yvec}_{(j)} \right)    \yvec_j^{\top} \mathcal{D}_{\mathcal{S}^{(12)}_{(i,j)}}^{-1}  \left( \xvec_i - \overline{\xvec}_{(i)} \right)    \right\}, \\
\end{align*}
are shown to be ratio consistent for the corresponding terms in ${\rm var}\left( \mathcal{U}_n \right)$. Standardizing $\mathcal{U}_n$ by the variance estimator, the Park-Ayyala test statistic is given by
\begin{equation}
T_{PA} = \frac{ \mathcal{U}_n}{\sqrt{ \left( \frac{n + m - 6}{n + m - 4} \right)^2 \left\{ \frac{2}{n(n - 1)} \widehat{\tr \left( \mathcal{R}_1^2 \right)} + \frac{2}{m(m-1)} \widehat{\tr \left( \mathcal{R}_2^2  \right) }+ \frac{4}{nm} \widehat{\tr \left( \mathcal{R}_1 \mathcal{R}_2  \right)}    \right\}  }}
\label{eqn:parkayyala}
\end{equation}

Asymptotic normality of the test statistic was derived under the following assumptions:
\begin{enumerate}[label=(PA \Roman*)]
	\addtocounter{enumi}{1}
	\item \label{itm:pa2} Sample sizes of the two groups, $n$ and $m$ are constrained as in \ref{itm:bs2}.
	\item \label{itm:pa3} If $\mathcal{R}$ is the correlation matrix, then $\tr \left( \mathcal{R}^4 \right) = o \left\{ \tr^2 \left(\mathcal{R}^2 \right) \right\}$. This condition is similar to \ref{itm:cq3}.
	\item \label{itm:pa4} The means $\muvec_1$ and $\muvec_2$ satisfy the local alternative condition $ n \left( \muvec_1 - \muvec_2 \right)^{\top} \mathcal{D}_{\mathcal{S}}^{-1/2} \mathcal{R} \mathcal{D}_{\mathcal{S}}^{-1/2} \left( \muvec_1 - \muvec_2 \right)^{\top} = o \left\{ \tr^2 \left( \mathcal{R}^2 \right) \right\}$
\end{enumerate}
The assumptions in \ref{itm:pa2} - \ref{itm:pa4} are milder than \ref{itm:sd1}-\ref{itm:sd4} and hold for a much larger family of covariance structures. Another major advantage of $T_{PA}$ is that it does not require the normality assumption. Instead, the test is constructed assuming the factor model defined in equation  \eqref{eqn:factormodel}.

The four test statistics have several key differences regarding their properties and performance. The Bai-Saranadasa test and Chen-Qin test are {\it orthogonal-transform} invariant, i.e. the operation $\xvec_i \mapsto \mathbf{U} \xvec_i, i = 1, \ldots, n$ and $\yvec_j \mapsto \mathbf{U} \yvec_j, j = 1, \ldots, m$ for some $p \times p$ orthogonal matrix $\mathbf{U}$ does not affect the test. The Srivastava-Du test and Park-Ayyala test are {\it scale-transform} invariant, wherein the operation described above does not affect the test when $\mathbf{U} = {\rm diag}\left(u_1, \ldots, u_p\right)$ is a diagonal matrix. In practice, scale transformation invariance is more useful than its orthogonal counterpart as they can bring variables on to a uniform scale. To better understand this difference, consider the contribution of each element towards the expected difference under the alternative when $\muvec_1 - \muvec_2 = \boldsymbol{\delta}$. In $T_{BS}$ and $T_{CQ}$, $k^{\rm th}$ element has a contribution of $\delta_k^2$, whereas in $T_{SD}$ and $T_{PA}$ the contribution is $\delta_k^2/\sigma_{kk}$. While the former depends on the scale of the variable, the latter is the coefficient of variation and is hence scale-free. In a scenario where the non-zero $\delta_k$'s correspond only to the values whose means are small, then $T_{PA}$ and $T_{SD}$ have higher power of detecting the difference. 

Due to their similarities in construction and assumptions, $T_{CQ}$ and $T_{PA}$ are observed to be applicable to a broader range of models. This is mainly because of relaxed assumptions on the covariance structure and lack of direct relationship between $p$ and $n$. However it is worth noting that the assumptions \ref{itm:bs1} and \ref{itm:sd1} are {\it asymptotic} and cannot be validated from a finite sample data set. For example, a data set with $p = 10,000$ and $n = 10$ can either imply the rate is polynomial ($p = n^4$) or linear ($p = 1000n$). There is no practical means of determining the true rate with a single data set. Another aspect of this asymptotic rate that is worth considering is that the number of variables is generally deterministic. In genomics data sets such as DNA methylation or gene expression, the dimension is the number of genes, which is fixed. The sample size is the number of biological replicates, which can be increased by collecting more specimens. Hence rate of increase cannot be used as a means to prefer one test to the other. A better approach to determine which method best suits a data set is through a simulation study. A controlled simulation study should be designed using the properties of the data set such as dependence structure and sparsity. The empirical type I error obtained by specifying equal means can be used to compare the performance of the methods. This approach was used in Ayyala \etal \citep{Ayyala2015} to determine that $T_{CQ}$ outperforms the other tests at controlling type I error rate and achieves reasonable power for immuno-precipitation based DNA methylation data. 

\subsection{Projection based tests} 
\label{sec:projectionmvt1}
The driving motivation behind the tests in Section \ref{sec:asymptoticmvt} is the fact that when $p > n$, the Hotelling's $T^2$ test statistic cannot be calculated. Another approach that has been considered to overcome this problem is to project the data into a lower dimensional space such that the assumptions of Hotelling's $T^2$ are satisfied. For $k < p$, consider a matrix $\mathcal{R} \in \mathbb{R}^{k \times p}$ with full row rank. The difference of means, $\muvec_1 - \muvec_2$, when projected onto the column space of $\mathcal{R}$, is equal to zero if and only if the difference itself is zero,
$\mathcal{R} \left( \muvec_1 - \muvec_2 \right) = \mathbf{0} \,\,\,  \Leftrightarrow  \,\,\, \muvec_1 - \muvec_2 = \mathbf{0}.$
By this equivalence, the hypothesis in \eqref{eqn:meanhyp1} is equivalent to  
\begin{equation}
H_{0:\mathcal{R}} : \mathcal{R} \muvec_1 = \mathcal{R} \muvec_2 \hspace{5mm} \mbox{vs. } \hspace{5mm} H_{A:\mathcal{R}} : \mathcal{R}\muvec_1 \neq \mathcal{R} \muvec_2.
\label{eqn:meanhyp2}
\end{equation} 
The equivalence holds for all rank-sufficient matrices and for all dimensions with $k \leq p$, which is an extremely large collection of matrices. In practice, we can only evaluate it for a very small set of matrices, based on which the conclusion can be drawn. Hence the two key factors that will affect the result of the test are $k$ and construction of $\mathcal{R}$.

A natural choice of $\mathcal{R}$ for dimension reduction is using principal component analysis. Let $\Sigma = \mathbf{V} \Omega \mathbf{V}^{\top}$  be the eigenvalue decomposition of the common covariance matrix $\Sigma$. The matrix $\mathbf{V} = \left( \mathbf{v}_1, \ldots, \mathbf{v}_p \right)$ is orthogonal where the columns are the eigenvectors and $\Omega = {\rm diag}(\omega_1, \ldots, \omega_p)$ is the diagonal with eigenvalues along the diagonal. Eigenvalue decomposition of the pooled sample covariance matrix $\mathcal{S}$ yields
\begin{equation*}
\mathcal{S} = \mathbf{U} \Lambda \mathbf{U}^{\top}, \hspace{5mm} \mathbf{U} = \left( \mathbf{u}_1, \ldots, \mathbf{u}_p \right), \quad \Lambda = {\rm diag}(\lambda_1, \ldots, \lambda_p),
\end{equation*}	
where $\lambda_1 \geq \cdots \geq \lambda_p$ are the eigenvalues and $\mathbf{u}_1, \ldots, \mathbf{u}_p$ are the orthogonal eigenvectors.  Properties of the eigenvalues and eigenvectors will be discussed in detail in later sections. The eigenvalues give a measure of the amount of variability in the data in the direction of the corresponding eigenvector. The cumulative relative variance of any set of eigenvectors $\{\mathbf{u}_{a_1}, \ldots, \mathbf{u}_{a_m} \}$ is given by $\left( \lambda_{a_1} + \ldots + \lambda_{a_m} \right)/ \left( \lambda_1 + \ldots + \lambda_p \right)$. Any set of $k$ eigenvectors can be used to construct a $k$-dimensional subspace to project the data.  However using the first $k$ eigenvectors is most informative as it contains the maximum cumulative relative variance in the data, equal to $(\lambda_1 + \ldots + \lambda_k)/(\lambda_1 + \ldots + \lambda_p)$. Define the matrix $\mathbf{U}_{(k)} = \left(\mathbf{u}_1, \ldots, \mathbf{u}_k \right)$ of dimension $p \times k$ using the first $k$ columns of $\mathbf{U}$ and the projections as
\begin{equation*}
\xvec_i \mapsto \xvec^{*}_i = \mathbf{U}_{(k)}^{\top} \xvec_i, \quad i = 1, \ldots, n,
 \hspace{1cm} \yvec_j \mapsto \yvec_j^{*} = \mathbf{U}_{(k)}^{\top} \yvec_j, \quad j = 1, \ldots, m.
\end{equation*}
The sample means of the projected observations will be $\overline{\xvec}^* = \mathbf{U}_{(k)}^{\top} \overline{\xvec}$ and $\overline{\yvec}^* = \mathbf{U}_{(k)}^{\top} \overline{\yvec}$ respectively. The pooled sample covariance matrix of $\xvec^*$ and $\yvec^*$ is $\mathbf{U}_{(k)}^{\top} \mathcal{S} \mathbf{U}_{(k)}$, which, by orthogonality of the columns of $\mathbf{U}_{(k)}$, is a diagonal matrix given by $\Lambda_{(k)} = {\rm diag}(\lambda_1, \ldots, \lambda_k)$. For any $k$, we can calculate the Hotelling's $T^2$ statistic using the projected data as
\begin{align}
T^2_{Hot(k)} & = \frac{n + m - k - 1}{(n + m - 2)k} \, \frac{n m}{n + m} \left(\overline{\xvec}^* - \overline{\yvec}^* \right)^{\top} \Lambda_{(k)} ^{-1} \left(\overline{\xvec}^* - \overline{\yvec}^* \right) \nonumber \\
& = \frac{n + m - k - 1}{(n + m - 2)k} \, \frac{n m}{n + m} \, \mathop{\sum}_{j = 1}^k \frac{ \left( \overline{X}^*_j - \overline{Y}^*_j \right)^2}{\lambda_j}. \label{eqn:hotelingt2proj}
\end{align}
When $k = p$ and $p < n + m - 2$, we have the original Hotelling's $T^2$ statistic as defined in \eqref{eqn:hotelingt2}, i.e.  $T^2_{Hot(p)} = T^2_{Hot}$. For any $k \leq p$, the null distribution of $T^2_{Hot(k)}$ will be $F_{k, n + m - k -1}$.  

\begin{figure}[!h]
	\centering
	\includegraphics[width=1.0\linewidth]{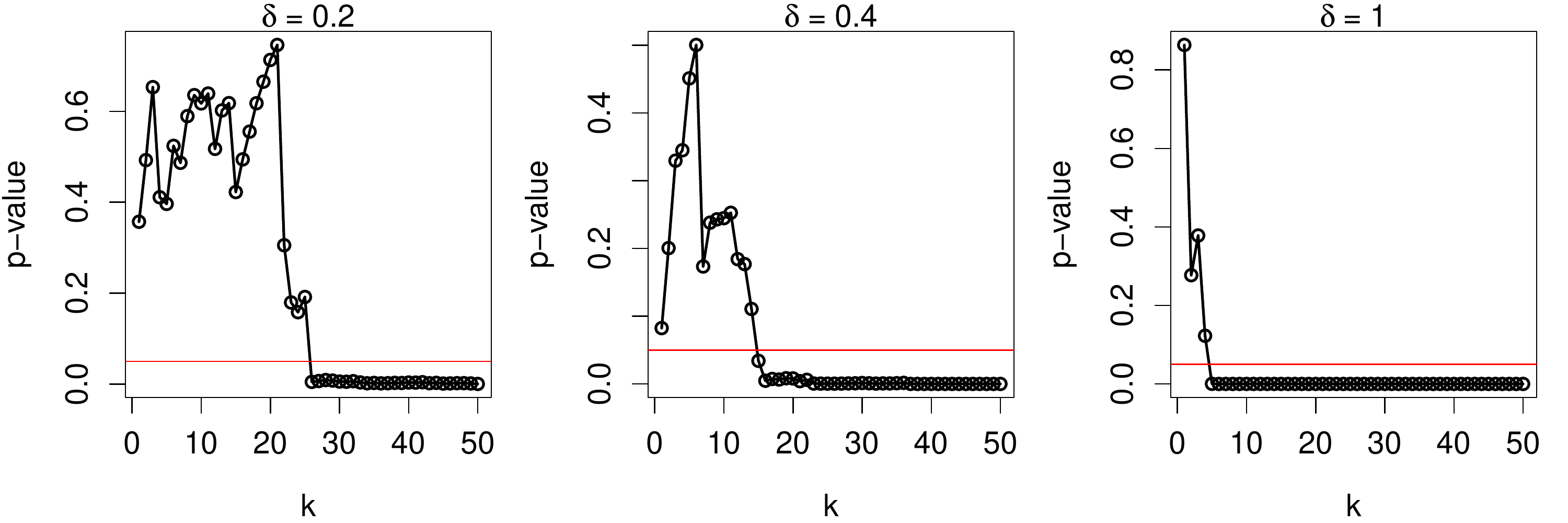}
	\caption{Figure}
	\label{fig:hotelingproj}
\end{figure}

While the motivation of projeting the samples into the principal component subspace is to reduce dimension and be able to use the Hotellings $T^2$ test statistic, one needs to be careful when choosing $k$. If the alternative hypothesis is true, then choosing a small $k$ can potentially lead to a type II error. This is because in $T^2_{Hot(k)}$, the summation will include only the first $k$ terms corresponding to the largest $\lambda$. But if the difference between the means is uniform over all the components, then the ratio of $ \left(\mu_{1j} - \mu_{2j} \right)/\lambda_j$ will be highlighted only for small $\lambda_j$, which correspond to large $j$. To illustrate this behaviour, consider the following study. Random samples are generated using $n = m = 100, p = 50$ and $\Sigma = {\rm diag}(\sigma_1, \ldots, \sigma_p)$ where $\sigma_i \sim {\rm Unif}(2,3)$. For the means, specify $\muvec_1 = (0, \ldots, 0)$ and $\muvec_{2} \sim (\delta, \ldots, \delta)$. Figure \ref{fig:hotelingproj} shows the $p$-value for different values of $\delta$ and for all $k \leq p$. 
The first thing to note is that $T^2_{Hot}$ detects the difference for the complete data ($k = p$), whereas projecting into a single dimension always fails to reject $H_0$. The smallest $k$ for which the p-value supports rejecting $H_0$ for $\delta = 0.2, 0.4$ and $1$ are $26, 15$ and $5$ respectively. The variance is kept constant for the three models, which implies the difference in results in due to $\delta$. As $\delta$ increases, there is greater separation between the two means and hence smaller $k$ is sufficient. The converse - rejecting the hypothesis for small $k < p$ when $\delta = 0$ and $H_0$ is rejected for $k = p$, is very unlikely to happen. Thus the type I error will be preserved for all $k$ but the projection test will have lower power than $T^2_{Hot}$.

In high dimensional setting, when $p > n + m - 2$, the eigenvalue matrix $\Lambda$ is singular with only the first $n + m - 2$ eigenvalues non-zero. Defining a {\it generalized} inverse of $\Lambda$ as ${\rm diag}(\lambda_1, \ldots, \lambda_{n + m - 2}, 0, \ldots, 0)^{-1} = {\rm diag}(1/\lambda_1, \ldots, 1/\lambda_{n + m - 2}, 0, \ldots, 0)$, the projected Hotelling's $T^2$ test statistic defined in \eqref{eqn:hotelingt2proj} can be calculated when $k \leq n + m - 2$. The full possible model, $T^2_{Hot(n + m - 2)}$ is not the complete Hotelling's $T^2$ test as it is still contains only a partial summation of terms in \eqref{eqn:hotelingt2proj}. However the $p$-value of $T^2{Hot(k)}$ behaves differently for different values of $k$ in high dimensions. In Figure \ref{fig:hotelingproj}, we observed that the type II error of $T^2_{Hot(k)}$ decreases as $k$ increases. This is because the deviations corresponding to the smallest eigenvalues will be included in the summation for large enough $k$, resulting in an increase in the test statistic value. But in high dimensions the smallest eigenvalues are set to zero. Therefore the projected $T^2$ test statistic can never achieve the value of $T^2_{Hot}$, resulting in extremely high type II error rate even for $k = n + m - 2$. 

\begin{figure}[!h]
	\centering
	\includegraphics[width=1.0\linewidth]{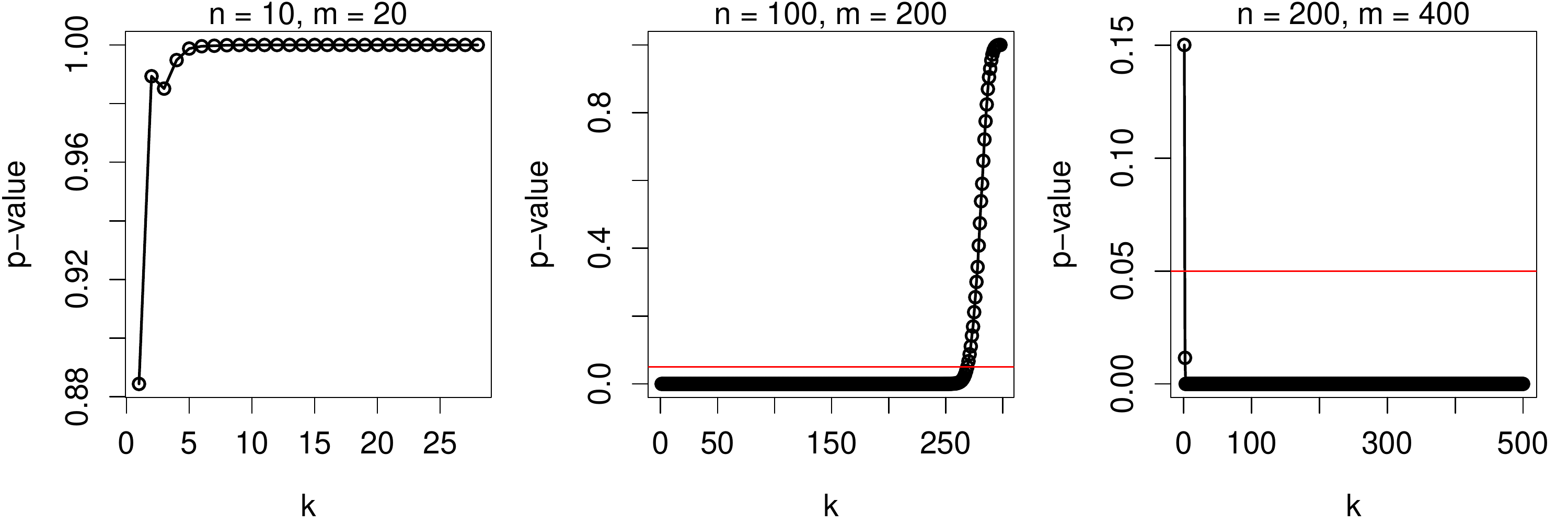}
	\caption{Figure}
	\label{fig:hotelingproj2}
\end{figure}

To illustrate the effect of $k$ on $T^2_{Hot(k)}$, consider the following simulation study. We set $p = 500$ and varied the sample sizes as $n \in \{10, 100, 200\}$ and $m = 2n$. The mean vectors are set as $\muvec_1 = (0, \ldots, 0)$ and $\muvec_2 = (1, \ldots, 1)$ respectively. 
The $p$-values of $T^2_{Hot(k)}$ for the three models and different values of $k$ within each model are presented in Figure \ref{fig:hotelingproj2}. Note that when $p > n + m$ as in the first two sub-figures, the $p$-value increases with $k$. Whereas in the third sub-figure, the properties of the $p$-value curve are similar to those seen in Figure \ref{fig:hotelingproj}. Similar results have been observed in a wide range of simulation models. Theoretical justification for this behaviour of the projection-based Hotelling's $T^2$ test is still lacking. 

\subsection{Random projections}
\label{sec:projectionmvt2}
As seen in the previous section, projecting on to the eigenspace of the pooled sample covariance matrix has its limitations in high dimension. The results also indicate that the conclusion will be contrary to the truth when sample sizes are small. While the concept of dimension reduction is effective, PCA is not the best approach for this task. Alternatively, projections based on {\it random matrices} have been shown to have good performance. A random projection embeds the $p$-dimensional variables into a lower $k$-dimensional space ($k << p$) while preserving the distances between points. The seminal result that allows us construct such an embedding is the Johnson-Lindenstrauss lemma \citep{Johnson1984}. 
\begin{theorem}[Johnson-Lindenstrauss lemma]
	For any collection of points $\xvec_1, \ldots, \xvec_n \in \mathbb{R}^p$ and $0 < \varepsilon < 1$, there exists a $k \geq k_0 = O \left(  \varepsilon^{-2} \log n \right)$ and a linear map $\mathcal{R}: \mathbb{R}^p \rightarrow \mathbb{R}^k$ such that
	\begin{equation}
	\left(1 - \varepsilon \right) \left\| \xvec_i - \xvec_j \right\|^2_2 \leq \left\| \mathcal{R}(\xvec_i) - \mathcal{R}(\xvec_j) \right\| \leq \left(1 + \varepsilon \right) \left\| \xvec_i - \xvec_j \right\|^2_2.
	\end{equation}
\end{theorem}
\noindent The theorem provides a method to determine the smallest dimension into which the original data can be embedded without altering the local properties of the data sets. Significance of this result is greatly enhanced by the fact that the dimension of the embedded space, $k$, depends on the sample size $n$ and not on the dimension $p$. While the result holds for any linear map, the most commonly used mapping is $\xvec \mapsto \mathcal{R} \xvec$ for some $k \times p$ matrix $\mathcal{R}$. To avoid the pitfalls of principal component based projections, the alternative is to randomly generate the matrix independent of the data.  

For a given $k \in \mathbb{Z}_{+}$, a random projection matrix $\mathcal{R} \in \mathbb{R}^{k \times p}$ is a matrix whose elements are random variables. Two distinctions need to be made when calling them {\it random} and {\it projection} matrices. Firstly, unlike matrices generated from distributions the matrix space such as Wishart, these random matrices are not structured. Secondly, these matrices need not necessarily have the properties of a projection matrix, {\it viz.} orthogonality, idempotency, etc. For simplicity of generation, the variables are chosen to be independent and identically distributed. Additional conditions can be imposed to provide structure to the projected data. While any distribution can be used to generate the elements, one property that is desired is that it is symmetric with zero mean and unit variance. This property ensures that the Euclidean distance between a pair of observations is preserved {\it in expectation}. That is, if $\mathbf{u} = (u_1, \ldots, u_p)$ is a random vector with $\mathbb{E}(u_i) = 0$ and $\mathbb{E}(u_i^2) = 1$, then 
\begin{equation*}
	\mathbb{E} \left( \|\mathbf{UX} - \mathbf{UY} \|_2^2 \right) = \mathbb{E} \left\{ \mathop{\sum}_{k = 1}^p u_k^2 \left(x_k - y_k \right)^2 \right\} = \mathop{\sum}_{k = 1}^p \mathbb{E} \left(u_k^2\right) \left(x_k - y_k \right)^2 = \| \mathbf{X} - \mathbf{Y} \|_2^2.
\end{equation*}
The most trivial distribution that is symmetric around zero with unit variance is the standard normal distribution, $\mathcal{N}(0, 1)$. To further simplify random number generation, one can also consider a uniform distribution ${\rm Unif} \left(-\sqrt{3}, \sqrt{3} \right)$, where the limits are adjusted to satisfy the moment conditions. 

Another class of projection matrices that has gained prominence recently is based on binary coins. Developed by Achlioptas \citep{Achlioptas2003}, this method is found to be very useful for dimension reduction in machine learning \citep{Fradkin2003}, image processing \citep{Bingham2001}and language processing \citep{Nunes2018}. The idea is to generate the elements of $\mathcal{R}$ from the set $\Omega = \{-1, 0, +1 \}$. Two distributions can be defined on $\Omega$ with zero mean and unit variance,
\begin{equation}
r_{ij}^{(1)} = \left\{  \begin{tabular}{lr} +1 & \mbox{ with probability   } 1/2 \\ 
												 -1  & \mbox{ with probability   }  1/2 \end{tabular} \right.,
\hspace{8mm}
r_{ij}^{(3)} = \left\{  \begin{tabular}{lr} $+\sqrt{3}$ & \mbox{ with probability   } \quad  1/6 \\ 
													0 	       & \mbox{ with probability   } \quad  2/3 \\
												 $-\sqrt{3}$  & \mbox{ with probability   } \quad  1/6 \end{tabular} \right..
\label{eqn:achlioptas}
\end{equation}
Among these two distributions, $r^{(3)}_{ij}$ is preferred to $r_{ij}^{(1)}$ because it produces a {\it sparse} embedding. By construction, the contribution of two out of three variables (on average) will be set to zero. Furthermore, the computation time is significantly improved when using $r_{ij}^{(3)}$. Extending from this work, Li \etal \citep{Li2006} generalized the procedure to define the distribution for any $\theta > 0$,
\begin{equation}
r_{ij}^{(\theta)} = \left\{  \begin{tabular}{ll} $+\sqrt{\theta}$ & \mbox{ with probability   } \quad  $\frac{1}{2 \theta}$  \\ 
															0 	       & \mbox{ with probability   } \quad  $1 -\frac{1}{\theta}$ \\
															$-\sqrt{\theta}$  & \mbox{ with probability   } \quad $ \frac{1}{2 \theta}$ \end{tabular} \right..
\label{eqn:pingli}
\end{equation}
This generalization improves on $r_{ij}^{(3)}$ as defined in \eqref{eqn:achlioptas} by increasing the sparsity of $\mathcal{R}$ with $\theta$, thereby reducing the computation cost of the projection. Li \etal \citep{Li2006} have shown that using $\theta$ as large as $p/\log(p)$ significantly reduces the computation cost with minimal loss of information (accuracy). However keeping in mind this trade-off between speed and information loss, the authors recommend $\theta = \sqrt{p}$.

Given a random projection matrix $\mathcal{P}$, the projected variables $\xvec^* = \mathcal{R} \xvec$ and $\yvec^* = \mathcal{P} \yvec$ have means $\mathcal{R} \muvec_1$ and $\mathcal{P} \muvec_2$ respectively. If the two populations are homogeneous, the common covariance matrix will be $\mathcal{R} \Sigma \mathcal{R}^{\top}$. Additionally if the variables are normally distributed, then the distribution is also preserved, i.e. $\xvec^* \sim \mathcal{N} \left( \mathcal{R} \muvec_1, \mathcal{R} \Sigma \mathcal{R}^{\top} \right)$ and $\yvec^* \sim \mathcal{N} \left( \mathcal{R} \muvec_2, \mathcal{R} \Sigma \mathcal{R}^{\top} \right)$. The sample means of the two populations will be $\overline{\xvec}^* = \mathcal{R} \overline{\xvec}$ and $\overline{\yvec}^* = \mathcal{R} \overline{\yvec}$ respectively and the pooled sample covariance matrix is $\mathcal{S}^* = \mathcal{R} \mathcal{S} \mathcal{R}^{\top}$. If $k < n + m - 2$, the Hotelling's $T^2$ test statistic for the projected data can be defined as 
\begin{align}
	T^2_{\mathcal{R}} & = \frac{n + m - k - 1}{(n + m - 2)k} \frac{ nm}{n + m} \left( \overline{\xvec}^* - \overline{\yvec}^* \right)^{\top} \left( \mathcal{R S R}^{\top}  \right)^{-1} \left( \overline{\xvec}^* - \overline{\yvec}^* \right) \nonumber \\
	& = \frac{n + m - k - 1}{(n + m - 2)k} \frac{ nm}{n + m} \left( \overline{\xvec} - \overline{\yvec} \right)^{\top}  \mathcal{R}^{\top} \left( \mathcal{R S R}^{\top}  \right)^{-1} \mathcal{R} \left( \overline{\xvec} - \overline{\yvec} \right).
\end{align}
Under the null hypothesis $H_{0:\mathcal{R}}$ defined in \eqref{eqn:meanhyp2}, $T^2_{\mathcal{R}}$ follows a $F_{k, n + m - k - 1}$ distribution conditional on $\mathcal{R}$. The p-value of the test statistic will be
\begin{equation}
	p_{\mathcal{R}} = 1 - F_{k, n + m - k - 1}(T^2_{\mathcal{R}})
	\label{eqn:projpvalue}
\end{equation}
At significance level $\alpha$, the null hypothesis is rejected if $p_{\mathcal{R}} < \alpha$.

In an unpublished work, Lopes \etal \citep{Lopes2015} first proposed \eqref{eqn:projpvalue} and suggested using $k = \lfloor (n + m)/2 \rfloor$ (assuming $p > \lfloor (n + m)/2 \rfloor$) for the dimension of the reduced space. They provide theoretical justification of conditions in which $T^2_{\mathcal{R}}$ has greater power than $T_{CQ}$ and $T_{SD}$. The only criticism of their procedure is the choice of $\mathcal{R}$. As the test statistic and $p$-value are calculated {\it conditional} on $\mathcal{R}$, the result of the test will be determined by the choice of the projection matrix $\mathcal{R}$. The results based on different realizations of the projection matrix $\mathcal{R}_1$ and $\mathcal{R}_2$ need not necessarily be consistent. To get rid of this sampling artefact, one should generate multiple instances of the projection matrix and combine the $p$-values of all the instances to draw inference. An exact method for combining the $p$-values from different projection matrices was developed by Srivastava \etal \citep{Srivastava2014}. Their method, RAPTT (stands for {\bf RA}ndom {\bf P}rojection {\bf T}-{\bf T}est)), uses the average $p$-value from multiple independent projection matrices to accept or reject the null hypothesis. The method works as follows. 

Consider $N$ random projection matrices $\mathcal{R}_1, \ldots, \mathcal{R}_N$ generated independently and the corresponding $p$-values calculated using equation \eqref{eqn:projpvalue},
\begin{equation*}
p_j = p_{\mathcal{R}_j} = 1 - F_{k, n + m - k - 1}(T^2_{\mathcal{R}_j}), \hspace{5mm} j = 1, \ldots, N.
\end{equation*}
Then the average $p$-value, $\overline{p} = N^{-1} \sum\limits_{j = 1}^N p_j$, is used to reject the null hypothesis at level $\alpha$. If $\psi_{\alpha}$ a cut-off based on the null distribution of $\overline{p}$ such that $P \left( \overline{p} > \psi_{\alpha} | H_0 \right)  = 1 - \alpha$, then we reject $H_0$ if $\overline{p} > \psi_{\alpha}$. The cut-off $\psi_{\alpha}$ is obtained from the null distribution of $\overline{p}$, which is not straightforward to derive. Instead, Srivastava \etal \citep{Srivastava2014} have established that the null distribution is independent of the parameters $\muvec_1, \muvec_2$ and $\Sigma$. Hence, without loss of generality, the values $\muvec_1 = \muvec_2 = \mathbf{0}$ and $\Sigma = \mathcal{I}_p$ can be used to derive the null distribution. Using this property, they proposed computing the cut-off empirically using the following algorithm.
\begin{enumerate}[label=(RAPTT \Roman*)]
	\item \label{itm:raptt1} Randomly generate $\xvec_1, \ldots, \xvec_n \sim \mathcal{N}\left(\mathbf{0}, \mathcal{I}\right)$ and $\yvec_1, \ldots, \yvec_m \sim \mathcal{N}\left(\mathbf{0}, \mathcal{I}\right)$.
	\item \label{itm:raptt2} Randomly generate $N$ projection matrices and calculate the $p$-values $p_1, \ldots, p_N$ using equation \eqref{eqn:projpvalue}. Calculate $\overline{p} = N^{-1} \sum\limits_{j = 1}^N p_j$. 
	\item \label{itm:raptt3} Repeat \ref{itm:raptt1} and \ref{itm:raptt2} $M$ times to calculate $\overline{p}_1, \ldots, \overline{p}_M$. Sort them in increasing order such that $\overline{p}_{[1]} \leq \ldots \leq \overline{p}_{[M]}$. Then the level $\alpha$ cut-off is estimated as
	\begin{equation}
	\widehat{\psi}_{\alpha} = \overline{p}_{\left[ M (1 - \alpha) \right]}.
	\end{equation}
\end{enumerate}

In their work, Srivastava \etal propose two types of projection matrices to use
\begin{enumerate}[label=(\roman*)]
	\item orthogonal matrices generated from the Haar distribution such that $\mathcal{R}\mathcal{R}^{\top} = \mathcal{I}_k$.
	\item a block-weighted approach where for each of the $k$ dimensions in the projected space, non-zero weights are assigned for a unique set of $[p/k]$ elements of the original variables. 
\end{enumerate}
A comprehensive simulation study reported in their work shows differences in the empirical power between the two projection matrices under certain situations. The type I error rates reported are relatively consistent. This discrepancy in the performance and its dependence on projection could be attributed to the limited scope of the simulation study (calculated based on 1000 runs). The optimal choice of projection matrix and a comprehensive investigation of performance of the projection-based tests under all models of projection matrix still needs to be addressed.

A major bottleneck of the projection-based tests is the computation time. Tests such as RAPTT are {\it exact} and are known to have better performance over asymptotic tests when the sample sizes are small. But the lack of null distribution and the variability of the test across different projection matrices imposes a heavy computational cost of the procedure. For instance, constructing the empirical null distribution using $N$ projection matrices and $M$ bootstrap samples for the data has a computational cost of $O(NM\tau)$, where $\tau$ is the  cost of calculating the Hotelling's $T^2$ test statistic and the corresponding $p$-value. Considering $N = M = 10^3$ leads to a cost of $O(10^6 \tau)$, which requires massively parallel computing to keep achieve reasonable computation time. Variability of $T^2_{R}$ and its $p$-value over the distribution of the projection matrices can be studied to determine the number of bootstrap samples required to achieve a specified level of accuracy in empirical calculations. 

\subsection{Other approaches}
In sections \ref{sec:asymptoticmvt}-\ref{sec:projectionmvt2}, the test statistics were based on the norm of the difference of mean vectors, either $(\muvec_1 - \muvec_2)^{\top} (\muvec_1 - \muvec_2)$ or $(\muvec_1 - \muvec_2)^{\top} \Sigma^{-1} (\muvec_1 - \muvec_2)$. Asymptotics and projection based methods are two major approaches commonly considered, but the hypothesis in \eqref{eqn:meanhyp1} can also viewed in a different light. Several tests have been proposed by either (i) aggregating the evidence across the individual elements or (ii) observing the maximum difference across the elements of $\xvec$ and $\yvec$. In this section, some aggregate tests based on univariate methods for individual elements are presented.

\begin{enumerate}
	\item \label{itm:pct} Pooled component test (PCT): Wu \etal \citep{Wu2006} proposed a test statistic which is applicable when there is missing data, i.e. one or more variables are not for all the observations. PCT requires the two groups to be homogeneous and normally distributed. The test statistic is obtained by averaging the squares of $t$-test statistics for the individual variables,
	\begin{equation}
	T_{PCT} = \frac{1}{p} \msum_{k = 1}^p \frac{ n_k m_k}{n_k + m_k} \frac{ \left( \overline{X}_{k} - \overline{Y}_k \right)^2}{S_{k}}
	\label{eqn:pct}
	\end{equation}
	where $n_k$ and $m_k$ are the number of observations for which the $k^{\rm th}$ variable is observed from the first and second samples respectively. The quantities $\overline{X}_k, \overline{Y}_k$ and $S_{k}$ are also similarly estimates of $\mu_{1k}, \mu_{2k}$ and $\Sigma_{kk}$ estimated using the observed values. Using the first two moments, the null distribution was established to be a scaled chi-squared, $T_{PCT} \stackrel{H_0}{\sim}c \chi^2_d$. The parameters $c$ and $d$ can be estimated from the individual $t$-test statistics to obtain the approximating null distribution. 
	\item \label{itm:gct} Generalized component test (GCT): Gregory \etal \citep{Gregory2015} proposed a test statistic for heterogeneous populations, replacing the pooled $t$-test statistic in $T_{PCT}$ with the unpooled test statistic,
	\begin{equation}
		T = \frac{1}{p} \msum_{k = 1}^p \frac{ \left( \overline{X}_k - \overline{Y}_k \right)^2}{ \frac{s^2_{Xk}}{n_k} + \frac{s^2_{Yk}}{m_k}},\hspace{1cm} 
		T_{GCT} = \frac{ \sqrt{p} \left\{ T - \left(1 + n^{-1} \widehat{a}_n + n^{-2} \widehat{b}_n \right) \right\} }{\widehat{\zeta}_n}
	\label{eqn:gct}
	\end{equation}
	where $\overline{X}_k, s_{Xk}, \overline{Y}_k, s_{Yk}$ are the mean and standard deviations of the $k^{\rm th}$ component for the two samples respectively and $n_k$ and $m_k$ are as defined in \ref{itm:pct}. The quantities $\widehat{a}_n$ and $\widehat{b}_n$ are obtained by combining the higher order moments of the elements of $\xvec$ and $\yvec$. The denominator $\widehat{\zeta}_n$ is estimated using a window-based aggregate of the autocovariance function across the elements. The test statistic is shown to be asymptotically normal. A key assumption of GCT is that the elements of $\xvec$ and $\yvec$ are ordered so that the autocovariance function across the elements diminishes with increasing lag (e.g. a moving average model). This assumption is very restrictive compared to the other tests. 
	\item Cai \etal \citep{Cai2014} developed a test based on the maximum scaled difference across the elements of the variables. Under the assumption that the populations are homogeneous and normally distributed, the test statistic is given by
	\begin{equation}
	T_{CLX} = \frac{ n m}{n + m} \mathop{\max}_{1 \leq k \leq p}  \frac{  \left\{ \widehat{\Omega} \left( \overline{\xvec} - \overline{\yvec} \right) \right\}_k^2}{ \frac{1}{n + m} \left[ n \mathcal{S}_{\xvec}(\widehat{\Omega})_{kk} + m \mathcal{S}_{\yvec} (\widehat{\Omega})_{kk} \right] },
	\end{equation}
	where $\mathcal{S}_{\xvec}(\widehat{\Omega})$ and $\mathcal{S}_{\yvec} (\widehat{\Omega})$ are the biased sample variance estimates of $\widehat{\Omega} \xvec_1, \ldots \widehat{\Omega} \xvec_n$ and $\widehat{\yvec}_1, \ldots, \widehat{\yvec}_m$ respectively. The precision matrix $\Omega = \Sigma^{-1}$ is estimated directly using constrained $\ell_1$-minimization for inverse matrix estimation (CLIME \citep{Cai2011}) to avoid inverting the singular pooled sample covariance matrix $\mathcal{S}$. Asymptotic null distribution of $T_{CLX}$ is shown to be an extreme value distribution of type I and a level $\alpha$ test rejects $H_0$ when $T_{CLX} \geq 2 \log p - \log \{ \log p\} - \log \pi - 2 \log \{ \log (1 - \alpha) \}$.
	\item Zoh \etal \citep{Zoh2018} have developed a Bayesian hypothesis for the hypothesis in \ref{eqn:meanhyp1} using Bayes factor. Under the assumption of homogeneous normal distributiosn for $\xvec$ and $\yvec$, they considered a Jeffrey's prior for $\left( \muvec, \Sigma \right)$ and a conjugate normal prior for $\boldsymbol{\delta} = \muvec_1 - \muvec_2$,
	\begin{equation*}
	\pi \left( \muvec, \Sigma \right) \propto \left| \Sigma \right|^{-(p+1)/2}, \hspace{1cm} \pi \left( \boldsymbol{\delta} | \Sigma, \tau_0 \right) \sim \mathcal{N} \left( \mathbf{0}, \tau_0 \Sigma \right), \tau_0 \in \mathbb{R}_+.
	\end{equation*}
	Then the Bayes factor was shown to admit a closed form given by
	\begin{equation}
	BF_{10} \left(\xvec, \yvec \right) = \frac{g \left( \xvec, \yvec | H_A \right)}{g \left( \xvec, \yvec | H_0 \right)} = \left(1 + \eta \right)^{-p/2} \left[ \frac{ 1 + \frac{p}{(1 + \eta)(n + m - p - 1)} T^2_{Hot}}{1 + \frac{p}{n + m - p - 1} T^2_{Hot}} \right]^{-(n + m - 1)/2}
	\label{eqn:zoh} 
	\end{equation}
	where $\eta = \frac{n m}{(n + m) \tau_0}$. They also proposed calculating the Bayes factor $BF_{10} \left( \mathcal{R} \xvec, \mathcal{R} \yvec \right)$ for any random projection matrix $\mathcal{R} \in \mathbb{R}^{k \times p}$ by replacing $p$ with $k$ and $T^2_{Hot}$ with $T^2_{Hot:\mathcal{R}}$ in \ref{eqn:zoh}. The rejection is constructed by translating the Hotelling's $T^2$ rejection region, $T^2_{Hot} > F_{\alpha, p, n + m - p - 1}$ to $BF_{10} \left(\xvec, \yvec \right)$. At significance level $\alpha$, the null hypothesis is rejected if 
	\begin{equation}
		BF_{10} \left(\xvec, \yvec \right) > \left( \tau_{\alpha}^* \right)^{-p/2} \left\{ 1 - \frac{\tau_{\alpha}^* - 1}{\tau_{\alpha}^*} C_n \right\},
	\end{equation}
	where $	C_n = (p F_{\alpha, p, n + m - p - 1}) \left\{ p F_{\alpha, p, n + m - p - 1} + n + m - p - 1 \right\}^{-1}$,
	$\tau_{\alpha}^* = n m \left\{ (n + m) \tau_{\alpha}  \right\}^{-1}$ and $\tau_{\alpha} = n m \left\{ (n + m) F_{\alpha, p, n + m - p - 1} - 1 \right\}^{-1}$.
\end{enumerate} 

\subsection{Dependent observations}
(write motivation)

For testing equality of means of two populations as presented in \eqref{eqn:meanhyp1}, the observations from each population are assumed to be independently and identically distributed. Most of the test statistics presented so far have been developed on several assumptions constraining the dependence structure. The testing problem has also been addressed when the covariance matrices are structured (Zhong \citep{Zhong2013}, Cai \citep{Cai2014}). But what happens if the observations are identically distributed but are not independent? Suppose the observations have the following covariance structure parametrized as $\cov\left( \xvec_i, \xvec_j \right) = \Sigma_1^{(i,j)}$ and $\cov \left( \yvec_i, \yvec_j \right) = \Sigma_2^{(i,j)}$. Then for any $i$ and $j$, the expected value of inner products of the variables will be $\mathbb{E} \left( \xvec_i^{\top} \xvec_j \right) = \muvec_1^{\top} \muvec_1 + \tr \left( \Sigma_1^{(i, j)} \right)$ and $\mathbb{E} \left( \yvec_i^{\top} \yvec_j \right) = \muvec_2^{\top} \muvec_2 + \tr \left( \Sigma_2^{(i, j)} \right)$ respectively. Considering the functional based on the Euclidean norm of $\overline{\xvec} - \overline{\yvec}$, its expected value will be
\begin{equation}
\mathbb{E} \left\{ \left( \overline{\xvec} - \overline{\yvec} \right)^{\top} \left( \overline{\xvec} - \overline{\yvec} \right)  \right\} = \left(\muvec_1 - \muvec_2 \right)^{\top} \left(\muvec_1 - \muvec_2 \right) + \frac{1}{n^2} \msum_{i, j = 1}^n \tr \left\{ \Sigma_1^{(i, j)} \right\} + \frac{1}{m^2} \msum_{i, j = 1}^m \tr \left\{ \Sigma_2^{(i, j)} \right\}.
\label{eqn:euclideandep1}
\end{equation}

Since the samples are assumed to be identically distributed, we have $\Sigma_1^{(i,i)} = \Sigma_1$, $\Sigma_2^{(i,i)} = \Sigma_2$. In the independent case, additionally we had $\Sigma_1^{(i, j)} = \Sigma_2^{(i, j)} = \mathbf{0}_{p \times p}$ when $i \neq j$. Under the dependence structure, we have additional $n(n - 1) + m(m - 1)$ covariance matrices in the model.  An unstructured dependence structure will therefore be infeasible because for any $i$ and $j$, we have only one pair of observations $(\xvec_i, \xvec_j)$ to estimate $\Sigma_1^{(i,j)}$. To make estimation feasible, we assume second-order stationarity on the dependence structures,
\begin{equation*}
\cov(\xvec_i, \xvec_j) = \Sigma_1^{(i,j)} = \Sigma_1 \left( i - j \right), \cov(\yvec_i, \yvec_j) = \Sigma_2^{(i,j)} = \Sigma_2 \left( i - j \right).
\end{equation*}
By symmetry, we have $\Sigma_1(-a) = \Sigma_1^{\top}(a)$ and $\Sigma_2(-a) = \Sigma_2^{\top}(a)$ for all $a \in \mathbb{Z}_{+}$. In time series, $\{\Sigma_1(a), a \in \mathbb{Z} \}$ and $\{\Sigma_2(a), a \in \mathbb{Z} \}$ represent the autocovariance functions of the two populations repsectively. The matrices $\Sigma_1(a)$ and $\Sigma_2(a)$ represent the autocovariance at lag $a$.

Using the autocovariance function, the expected value in \eqref{eqn:euclideandep1} simplifies to
\begin{align}
\mathbb{E} \left\{ \left( \overline{\xvec} - \overline{\yvec} \right)^{\top} \left( \overline{\xvec} - \overline{\yvec} \right)  \right\} & = \left(\muvec_1 - \muvec_2 \right)^{\top} \left(\muvec_1 - \muvec_2 \right) \nonumber \\
& + \frac{1}{n^2} \msum_{a = -(n-1)}^{n - 1}  \left(n - |a| \right) \tr \left\{ \Sigma_1(a)\right\} + \frac{1}{m^2} \msum_{a = -(m-1)}^{m - 1} \left( m - |a| \right)\tr \left\{ \Sigma_2(a) \right\}.
\label{eqn:euclideandep2}
\end{align}
A functional that is unbiased for the Euclidean norm of $\muvec_1 - \muvec_2$ can be constructed using \eqref{eqn:euclideandep2} as
\begin{equation}
\mathcal{M} =  \left( \overline{\xvec} - \overline{\yvec} \right)^{\top} \left( \overline{\xvec} - \overline{\yvec} \right)  - \left[ \frac{1}{n^2} \msum_{a = -(n-1)}^{n - 1}  \left(n - |a| \right) \tr \left\{ \widehat{\Sigma_1(a)}\right\} + \frac{1}{m^2} \msum_{a = -(m-1)}^{m - 1} \left( m - |a| \right)\tr \left\{ \widehat{\Sigma_2(a)} \right\} \right],
\end{equation}
where $\widehat{\Sigma_1(a)}$ and $\widehat{\Sigma_2(a)}$ are the biased estimators of $\Sigma_1(a)$ and $\Sigma_2(a)$ defined as
\begin{equation}
\widehat{\Sigma}_1(a) = \frac{1}{n} \msum_{i = 1}^{n - |a|} \left( \xvec_i - \overline{\xvec} \right) \left( \xvec_{i + a} - \overline{\xvec} \right)^{\top}, \hspace{5mm}
\widehat{\Sigma}_2(a) = \frac{1}{m} \msum_{i = 1}^{m - |a|} \left( \yvec_i - \overline{\yvec} \right) \left( \yvec_{i + a} - \overline{\yvec} \right)^{\top}.
\label{eqn:autocovestim}
\end{equation}
These estimators are the {\it biased} estimators (Brockwell and Davis \citep{Brockwell1986}), which should be of no concern to us since we are only interested in their trace.  When $p$ is finite, these estimators are known to be asymptotically unbiased. However in high dimensions, when $p$ increases with $n$ this property is no longer valid. For instance, the expected value of $\tr \left\{ \widehat{\Sigma}_1(a) \right\}$ will be  $\mathbb{E} \left[ \tr \left\{ \widehat{\Sigma}_1(a) \right\} \right] = \msum_{b = 0}^{n - 1} \theta_n(a,b) \tr \left\{ \Sigma_1(b) \right\}$, where 
\begin{align} 
\theta_n(a, b) & = \left(1 - \frac{a - 1}{n}\right) \mathbb{I}(a = b)+\left(1-\frac{a-1}{n}\right)\left(1-\frac{b-1}{n}\right) \frac{\{2-\mathbb{I}(a = 1)\}}{n} \nonumber \\ &-\frac{1}{n^{2}} \sum_{t=1}^{n-a+1} \sum_{s=1}^{n}\{\mathbb{I}(|t-s|+1 = b) + \mathbb{I}(|t+i-s-1| = b)\} 
\end{align} 

Asymptotic unbiasedness for finite $p$ follows from the leading term converging to 1 and the second and third terms, which are $O(n^{-1})$, converging to zero as $n$ goes to infinity because $\tr \left\{ \Sigma_1(a) \right\} = O(1)$. In high dimension, if the autocovariance structure is {\it proper} with all eigenvalues being non-zero, then $\tr \left\{ \Sigma_k(a) \right\} = O(p)$ for $k = 1, 2$ and all lags $a$. Hence all three terms in the expression for $\theta_n(a,b)$ should be considered. The expected value of $\tr \left\{ \widehat{\Sigma}_1(a) \right\}$ depends on the autocovariance matrices at all lags through the trace function, which is a univariate measure of the matrix. Expressing in vector form, we have $\mathbb{E} \left\{ \widehat{\boldsymbol{\gamma}_n} \right\} = \Theta_n \boldsymbol{\gamma}$ where $\Theta_n = (\theta_n(a,b))_{a, b \in \{0, \ldots, n-1\}}, \boldsymbol{\gamma} = \left( \tr \left\{ \Sigma_1(0) \right\}, \ldots, \tr \left\{ \Sigma_1(n - 1) \right\} \right)$ and $\widehat{\boldsymbol{\gamma}} = \left( \tr \left\{ \widehat{\Sigma}_1(0) \right\}, \ldots, \tr \left\{ \widehat{\Sigma}_1(n - 1) \right\} \right)$ respectively. This property can be used to construct unbiased estimators for $\tr \left\{ \Sigma_1(0) \right\}$ as elements of the vector $\widehat{\boldsymbol{\gamma}^*} = \Theta_n^{-1} \widehat{\boldsymbol{\gamma}}_n$. Denoting the elements of $\widehat{\boldsymbol{\gamma}^*}$ as $\widehat{\tr \left\{ \Gamma(a) \right\}}$, the functional can finally be constructed as
\begin{equation}
\mathcal{M}_n=  \left( \overline{\xvec} - \overline{\yvec} \right)^{\top} \left( \overline{\xvec} - \overline{\yvec} \right)  - \left[ \frac{1}{n^2} \msum_{a = -(n-1)}^{n - 1}  \left(n - |a| \right) \widehat{ \tr \left\{ \Sigma_1(a) \right\}} + \frac{1}{m^2} \msum_{a = -(m-1)}^{m - 1} \left( m - |a| \right) \widehat{\tr \left\{ \Sigma_2(a) \right\}} \right],
\label{eqn:dependmn}
\end{equation}

Ayyala \etal \citep{Ayyala2017} proposed a test statistic based on $\mathcal{M}_n$ defined in \eqref{eqn:dependmn}. In addition to the second-order stationary autocovariance structure, observations from the two populations are assumed to be realizations of two independent $M$-dependent strictly stationary Gaussian processes with means $\muvec_1$ and $\muvec_2$ and autocovariance structures $\{\Sigma_1(a)\}$ and $\{\Sigma_2(a)\}$ respectively. The $M$-dependence structures imposes the autocovariance matrices to be equal to zero for lags greater than $M$. Properties of the test statistic are established based on the following assumptions:
\begin{enumerate}[label=(APR \Roman*)]
	\item \label{itm:apr1} The observations are realizations of $M$-dependent strictly stationary Gaussian processes.
	\item \label{itm:apr2} The rates of increase of dimension $p$ and order $M$ with respect to $n$ are linear and polynomial respectively, 
	\begin{equation*}
		p = O(n), \hspace{1cm} M = O(n^{1/8}).
	\end{equation*}
	\item \label{itm:apr3} For any $k_1, k_2, k_3, k_4 \in \{1, 2\}$, 
	\begin{equation*}
	\tr \left\{ \Sigma_{k_1}(a) \Sigma_{k_2}(b) \Sigma_{k_3}(c) \Sigma_{k_3}(d) \right\} = o\left\{ (M + 1)^{-4} \tr^{2} \left(\Omega_{1}+\Omega_{2} \right)^{2}\right\},
	\end{equation*}
	where $\Omega_1 = \msum_{a = -M}^M (1 - |a|/n) \Sigma_1(a)$ and $\Omega_2 = \msum_{a = -M}^M (1 - |a|/n) \Sigma_2(a)$. 
	\item \label{itm:apr4} The means $\muvec_1$ and $\muvec_2$ satisfy the local alternative condition
	\begin{equation*}
		\left(\muvec_1 - \muvec_2 \right)^{\top} \left\{\Sigma_w(a) \Sigma_w(-a) \right\}^{\frac{1}{2}} \left(\muvec_1 - \muvec_2 \right) = o \left\{ (M + 1)^{-4} n^{-1} \tr \left( \Omega_1 + \Omega_2 \right)^{2} \right\}	
	\end{equation*}
\end{enumerate}

Setting $M = 0$ and $\Sigma_1(a) = \Sigma_2(a) = 0$ for all $a \neq 0$, it is straightforward to see that the conditions \ref{itm:apr3} and \ref{itm:apr4} are similar to \ref{itm:cq3} and \ref{itm:cq4}. The test statistic is given by
\begin{equation}
T_{APR} = \frac{\mathcal{M}_n}{\sqrt{ \widehat{ {\rm var} \left(\mathcal{M}_n \right)}}}
	\label{eqn:ayyalapr}
\end{equation}
where the variance estimate is constructed similar to $T_{CQ}$ and $T_{PA}$ using a leave-out method for better asymptotic properties. For exact form of the estimator, please refer to Ayyala \etal \citep{Ayyala2017}. Under the conditions \ref{itm:apr1}-\ref{itm:apr4}, $T_{APR}$ is shown to be asymptotically normal. While the test statistic and the empirical studies of Ayyala \etal are valid, Cho \etal \citep{Cho2019} identified some theoretical errors in the proofs and provided some corrections to some results and assumptions in Ayyala \etal. 

One issue that still needs to be addressed is the choice of $M$. Simulation studies reported in Ayyala \etal indicate that over-estimating $M$ is better than underestimating. When the specified value of $M$ in the analysis is greater than the true order of dependency, the error is in estimating zero matrices for lags greater than the true $M$. Under-specifying the value results in bias as autocovariances for several lags will not be estimated. Accurate estimation of $M$ using the data is not addressed and remains an open area of research. A large class of models can be approximated using $M$-dependent strictly stationary processes. Tests for other classes of models such as second-order stationary processes or Non-Gaussian processes is another area of active research. 

\section{Covariance matrix}
\label{sec:covariancematrix}
The covariance matrix of a multivariate random variable is a measure of dependence between the components of the variable. It is the second order central moment of the variable, defined as $\Sigma = {\rm var}(\xvec) = \mathbb{E} \left\{ \left(\xvec - \muvec \right) \left(\xvec - \muvec \right)^{\top} \right\}$, where $\muvec = \mathbb{E}(\xvec)$. The covariance matrix is often re-parameterized using its inverse, called the {\bf precision} matrix, $\Omega = \Sigma^{-1}$. Elements of the precision matrix are useful in determining conditional independence under normality. If $\xvec \sim \mathcal{N}(\muvec, \Sigma)$, then $\Omega_{ij} = 0$ implies $X_i$ is independent of $X_j$ conditional on $\{X_k : k \neq i, j\}$. The precision matrix is important because it can be used to construct an undirected graphical network model. Representing the components as nodes of the network, edges are defined by the elements of $\Omega = (\omega_{ij})$, where $\omega_{ij} \neq 0$ indicates the presence of an edge and $\omega_{ij} = 0$ indicates the absence of an edge between nodes $i$ and $j$. In view of these properties of the covariance matrix and other distributional properties, normality of the variables is commonly used in covariance matrix estimation. Unless otherwise stated, we shall assume the variables are normally distributed for the remainder of this section.  

Given an \iid sample $\xvec_i \sim \mathcal{N}(\muvec, \Sigma), i = 1, \ldots, n$, the biased sample covariance matrix is defined as
\begin{equation}
\mathcal{S} = \left(s_{ij} \right)_{i, j = 1, \ldots, p} = \frac{1}{n} \sum \limits_{i = 1}^n \left( \xvec_i - \overline{\xvec} \right) \left( \xvec_i - \overline{\xvec} \right)^{\top}
\label{eqn:samplecov}
\end{equation}
with $\mathbb{E}(\mathcal{S}) = (n - 1)/n \Sigma$ and ${\rm rank}(\mathcal{S}) = \min(n - 1, p)$. In traditional multivariate setting with $p < n$, $\mathcal{S}$ is non-singular and consistent for $\Sigma$. The sampling distribution of $\mathcal{S}$ is a Wishart distribution with $n - 1$ degrees of freedom (Anderson \citep{AndersonBook}, Muirhead \citep{MuirheadBook}). Additionally, the eigenvalues of $\mathcal{S}$ are also consistent for the eigenvalues of $\Sigma$. Asymptotically, the eigenvalues are normally distributed - a result that can be used to construct hypothesis tests. Estimation of eigenvalues of $\Sigma$ is of importance because they give the variance of the {\it principal components}, which are useful in constructing lower-dimensional embeddings of the data (dimension reduction). Hypothesis tests concerning the structure of the covariance matrix such as sphericity ($H_0 : \Sigma = \sigma^2 \mathcal{I}$) and uniform correlation ($H_0 : \Sigma = \sigma^2 \left[ (1 - \rho) \mathcal{I} + \rho \mathbf{1} \mathbf{1}^{\top} \right]$) are constructed using this property (\citep{AndersonBook, MuirheadBook}). Testing equality of covariance matrices for two or more groups is also well-defined when using the sample covariance matrix and its Wishart properties. 

Results from traditional multivariate analysis are valid only when $n > p$ and $p$ is assumed to be fixed. In high dimensional analysis, as seen in Section \ref{sec:meanvt}, $p$ is assumed to be increasing with $n$. How can we construct consistent estimators for $\Sigma$ and test statistics to compare the covariance structures of two or more populations in high dimension? In high dimensional models with $p \geq n$, the sample covariance matrix $\mathcal{S}$ is rank-deficient. Estimation of $\Sigma$ and $\Omega$ also suffer from the curse of dimensionality even when $p < n$ with $p/n \rightarrow c \in (0,1)$. When $p \rightarrow \infty$, $\mathcal{S}$ is no longer consistent for $\Sigma$. Estimation of $\Sigma$ was not an issue in tests the mean vector since we were only interested in consistent estimator for a function of $\Sigma$, e.g. $\tr \left(\Sigma \right)$ or $\tr \left( \Sigma^2 \right)$. 

\subsection{Estimation}

To obtain consistent estimators for $\Sigma$, two methods for reducing the parameter space dimension are used - structural constraints or regularization through sparsity. A banding approach, proposed by Bickel and Levina \citep{Bickel2008} sets elements outside a {\it band} around the diagonal to zero. For any $1 \leq k \leq p$, the banded estimator $\widehat{\Sigma}^{(k)}$ is defined as
\begin{equation}
\widehat{\Sigma}^{(k)}_{ij} = \left\{ \begin{matrix} s_{ij} & \mbox{ if } |i - j| < k \\ 0 & \mbox{ if } |i - j| \geq k \end{matrix}\right..
\end{equation}
Here $k$ denotes the width of the band, clearly indicating $\widehat{\Sigma}^{(k)}$ as the diagonal estimator. The estimator is consistent for $\Sigma$ under the $\ell_2$ matrix norm and when $\log p/n \rightarrow 0$. The optimal value of $k$ is chosen using $K$-fold cross validation of the estimated risk. It is particularly effective when the components of $\xvec$ are ordered so that $\sigma_{ij}$ decreases as $|i - j|$ increases. Consistency of the estimator is also shown to hold for non-Gaussian variables whose elements have sub-exponential tails

Regularization is a more commonly used approach for covariance matrix estimation as it is easier to formulate mathematically. Under normality, likelihood of $\Sigma$ given a sample $\xvec_1, \ldots, \xvec_n$ can be expressed as
\begin{align}
\mathcal{L} \left( \Sigma | \xvec_1, \ldots, \xvec_n \right) & = \msum \limits_{i = 1}^n - \log \left\{ \sqrt{{\rm det} \, \Sigma} \right\} - \left(\xvec_i - \overline{\xvec} \right)^{\top} \Sigma^{-1} \left( \xvec_i - \overline{\xvec} \right) \nonumber \\
& = -\frac{n}{2} \left[ \log \left\{ {\rm det} \, \Sigma \right\} + 2 \tr \left\{ \mathcal{S} \Sigma^{-1} \right\} \right].
\label{eqn:covariancelkhd}
\end{align}
Expression of the second term follows by applying the matrix result that for any $p$ dimensional vector $\mathbf{x}$ and $p \times p$ matrix $B$, we have $\mathbf{x}^{\top} B \mathbf{x} = \tr \left(  \mathbf{x}^{\top} B \mathbf{x}  \right) = \tr \left(  B \mathbf{x} \mathbf{x}^{\top} \right)$. Alternatively, the likelihood can be expressed in terms of the precision matrix $\Omega$ as 
\begin{equation}
\mathcal{L} \left( \Omega | \xvec_1, \ldots, \xvec_n \right) = \frac{n}{2} \left[ \log \left\{ {\rm det} \, \Omega \right\} - 2 \tr \left\{ \mathcal{S} \Omega \right\} \right].
\label{eqn:precisionlkhd}
\end{equation}
Maximizing the likelihood in \eqref{eqn:covariancelkhd} with respect to $\Sigma$ yields $\widehat{\Sigma} = \mathcal{S}$. 

Regularization of the covariance matrix estimator is achieved by adding a penalty term to the likelihood in \ref{eqn:covariancelkhd},
\begin{equation}
\mathcal{L}^*(\Sigma | \xvec_1, \ldots, \xvec_n) = -\frac{n}{2} \left[ \log \left\{ {\rm det} \, \Sigma \right\} + 2 \tr \left\{ \mathcal{S} \Sigma^{-1} \right\} \right] - \lambda \, \mathcal{P}(\Sigma), 
\label{eqn:covariancepenal1}
\end{equation}
for some penalty function $\mathcal{P}$ which can be defined to achieve a desired effect on $\widehat{\Sigma}$. The penalty parameter $\lambda$ dictates the trade-off between maximizing the likelihood term and minimizing the penalty. Inspired by lasso (Tibshirani \citep{Tibshirani1996}), Bien and Tibshirani \citep{Bien2011} proposed using a $\ell_1$-penalty to induce sparsity in the estimator. The penalty function is given by $\mathcal{P}(\Sigma) = \| \mathcal{W} \circ \Sigma \|_1 = \msum_{i,j} w_{ij} \sigma_{ij}$, where $\circ$ denotes the Hadamard element-wise product. The matrix $\mathcal{W} = \mathbf{1} \mathbf{1}^{\top}$ penalizes all the elements of $\Sigma$ whereas $\mathcal{W} = \mathbf{1} \mathbf{1}^{\top} - \mathcal{I}$ penalizes only the off-diagonal terms. Another approach to address regularization was developed by Daniels and Kass \citep{Daniels2001} by shrinking the eigenvalues to make the estimator more stable. 

While theoretically developing penalized estimates for the covariance matrix is important, it is practically more conducive to obtain sparse estimates of the precision matrix. Sparsity of precision matrix translates to absence of edges between nodes in the network model. Hence a sparse precision matrix can be used to isolate clusters of nodes which are strongly dependent within themselves and independent of the other clusters. The $\ell_1$ penalized precision matrix estimation is done by maximizing the function
\begin{equation}
\mathcal{L} \left( \Omega | \xvec_1, \ldots, \xvec_n \right) = \frac{n}{2} \left[ \log \left\{ {\rm det} \, \Omega \right\} - 2 \tr \left\{ \mathcal{S} \Omega \right\} \right] - \lambda \, \| \Omega \|_1.
\label{eqn:covariancepenal2}
\end{equation}
Termed by Friedman \etal \citep{Friedman2008} as {\bf glasso} (short for graphical lasso), the problem has garnered great levels of interest. Several extensions and improvisations of the original glasso method have been proposed. Danaher \etal \citep{Danaher2014} and Guo \etal \citep{Guo2011} studied joint estimation of $K > 1$ precision matrices by imposing two levels of penalties. For sparse estimation of precision matrices $\Omega^{(1)}, \ldots, \Omega^{(K)}$, using \eqref{eqn:covariancepenal2} individually will not preserve the cluster structure across the groups. By introducing a penalty to merge the $K$ groups, the following penalty functions have been proposed:
\begin{align*}
\mbox{Fused graphical lasso:} \hspace{5mm} & \mathcal{P}\left( \Omega^{(1)}, \ldots, \Omega^{(K)} \right) = \lambda_1 \sum \limits_{k = 1}^K \sum \limits_{i \neq j} |\omega^{(k)}_{ij}| + \lambda_2 \sum \limits_{k < m} \sum \limits_{i \neq j} | \omega^{(k)}_{ij} - \omega^{(m)}_{ij}|, \\
\mbox{Group graphical lasso:} \hspace{5mm} & \mathcal{P}\left( \Omega^{(1)}, \ldots, \Omega^{(K)} \right) = \lambda_1 \sum \limits_{k = 1}^K \sum \limits_{i \neq j} |\omega^{(k)}_{ij}| + \lambda_2  \sum \limits_{i \neq j} \left( \sum \limits_{k = 1}^K \omega^{(k)^2}_{ij}\right)^{1/2}, \\
\mbox{Guo \etal:} \hspace{5mm} & \mathcal{P}\left( \Omega^{(1)}, \ldots, \Omega^{(K)} \right) = \lambda_1 \sum \limits_{i \neq j} |\theta_{ij}| + \lambda_2  \sum \limits_{i \neq j} \sum \limits_{k = 1}^K |\gamma^{(k)}_{ij}|,
\end{align*}
where Guo \etal parameterize the $K$ precision matrices as $\Omega^{(k)} =\Theta \circ \Gamma^{(k)}$ with $\Theta$ representing the overall network structure and $\Gamma^{(k)}$'s representing the group-specific difference in the structure.

\subsection{Hypothesis testing}
When studying the covariance matrix of a multivariate Gaussian population, there are two common hypotheses of interest: 
\begin{align}
\mbox{Sphericity: } \hspace{1cm} & H_0 : \Sigma = \sigma^2 \mathcal{I} \hspace{5mm} \mbox{ vs. } \hspace{5mm} H_A : \Sigma \neq \sigma^2 \mathcal{I} \nonumber \\
\mbox{Identity: } \hspace{1cm} & H_0 : \Sigma = \mathcal{I} \hspace{5mm} \mbox{ vs. } \hspace{5mm} H_A : \Sigma \neq \mathcal{I}
\label{eqn:covmathyps}
\end{align}
These hypotheses can be alternatively stated using eigenvalues. If $\lambda_1,\ldots, \lambda_p$ are the eigenvalues of $\Sigma$, the hypotheses in \eqref{eqn:covmathyps} are equivalent to 
\begin{align}
\mbox{Sphericity: } \hspace{1cm} & H_0 : \lambda_1 = \ldots = \lambda_p \hspace{5mm} \mbox{ vs. } \hspace{5mm} H_A: \lambda_i \neq \lambda_j \mbox{ for some } i \neq j \nonumber \\
\mbox{Identity: } \hspace{1cm} & H_0 : \lambda_i = 1 \,\, \forall \,\, i = 1, \ldots, p \hspace{5mm} \mbox{ vs. } \hspace{5mm} H_A : \lambda_i \neq 1 \mbox{ for some } i.
\label{eqn:covmathyps2}
\end{align}
Functionals of $\Lambda = (\lambda_1, \ldots, \lambda_p)$ which are equal to zero under the null hypothesis can be constructed by observing that under sphericity, the {\it variance} of $\Lambda$ is equal to zero. For the identity hypothesis, deviation of $\Lambda$ from one is zero. The functionals (John \citep{John1972}, Nagao \cite{Nagao1973}) can be defined as 
\begin{gather}
U(\Lambda) = \frac{1}{p} \sum \limits_{k = 1}^p \left( \frac{\lambda_k}{\overline{\lambda}} - 1 \right)^2 = \frac{1}{p} \tr \left\{ \frac{\Sigma}{\tr\Sigma/p} - \mathcal{I}  \right\}^2, \nonumber \\
V(\Lambda) = \frac{1}{p} \sum \limits_{k = 1}^p \left( \lambda_k - 1\right)^2 = \frac{1}{p} \tr \left\{  \Sigma - \mathcal{I}  \right\}^2,
\label{eqn:covmatfunctionals}
\end{gather}

In the traditional setting when $p < n$, the sample covariance matrix $\mathcal{S}$ (and its eigenvalues) are consistent for $\Sigma$ (and $\Lambda$). Hence the test statistics based on functionals in \eqref{eqn:covmatfunctionals} are
\begin{equation}
U_n =  \frac{1}{p} \tr \left\{ \frac{\mathcal{S}}{\tr\mathcal{S}/p} - \mathcal{I}  \right\}^2, \hspace{1cm} V_n = \frac{1}{p} \tr \left\{  \mathcal{S} - \mathcal{I}  \right\}^2
\label{eqn:covmatteststat1}
\end{equation}
which are shown to follow chi-squared distributions asymptotically with $p(p+1)/2 - 1$ and $p(p+1)/2$ degrees of freedom respectively. Ledoit and Wolf \citep{Ledoit2002} studied the properties of $U_n$ and $V_n$ when $p/n \rightarrow c > 0$ and observed that $U_n$ performs well even in the high-dimensional case. For the identity hypothesis, they constructed a new test statistic,
\begin{equation}
W_n = \frac{1}{p} \tr \left\{  \mathcal{S} - \mathcal{I}  \right\}^2 - \frac{p}{n} \left[ \frac{1}{p} \tr \mathcal{S} \right]^2 + \frac{p}{n}
\label{eqn:covmatteststat2}
\end{equation}
which is also asymptotically chi-squared with $p(p+1)/2$ degrees of freedom but has better properties than $V_n$. Relaxing the assumption of normal distribution and a direct relationship between $n$ and $p$, Chen \etal proposed test statistics $U_n^*$ and $V_n^*$ which are asymptotically normally distributed. These test statistics are in the same spirit as $T_{CQ}$ \eqref{eqn:chenqin} and uses leave-out cross-validation type products to improve the asymptotic properties. 

Next, consider testing equality of covariance matrices from two normal populations $\xvec_i \sim \mathcal{N}(0, \Sigma_1), i = 1, \ldots, n$ and $\yvec_j \sim \mathcal{N}(0, \Sigma_2), j = 1, \ldots, m$. The sample covariance matrices and pooled covariance matrix,
\begin{equation*}
\mathcal{S}_1 = \frac{1}{n} \sum \limits_{i = 1}^n \left( \xvec_i - \overline{\xvec} \right) \left( \xvec_i - \overline{\xvec} \right)^{\top}, 
\mathcal{S}_2 = \frac{1}{m} \sum \limits_{j = 1}^m \left( \yvec_j - \overline{\yvec} \right) \left( \yvec_j - \overline{\yvec} \right)^{\top}, 
\mathcal{S}_{pl} = \frac{ n \mathcal{S}_1 + m \mathcal{S}_2}{n + m}
\end{equation*}
are used to construct the likelihood ratio test statistic as
\begin{equation}
\mathcal{L} = - \left\{ (n + m) \log |\mathcal{S}_{pl}| - n \log |\mathcal{S}_1| - m \log | \mathcal{S}_2|  \right\}
\end{equation}
Under $H_0 : \Sigma_1 = \Sigma_2$, $\mathcal{L}$ asymptotically follows a chi-squared distribution with $p(p + 1)/2$ degrees of freedom. Extending to $K$ groups, the test statistic is 
\begin{equation*}
\mathcal{L}_K = -\left\{ \sum \limits_{g = 1}^K n_g \left( \log | \mathcal{S}_{pl}| - \log | \mathcal{S}_g| \right)  \right\},
\end{equation*}
where $n_g$ is the sample size of the $g^{\rm th}$ group and $\mathcal{S}_{pl} = \left( \sum_{g=1}^K n_g \right)^{-1} \left(\sum_{g=1}^K n_g \mathcal{S}_g \right)$. Under $H_0 : \Sigma_1 = \ldots = \Sigma_K$, the LRT statistic $\mathcal{L}_K$ asymptotically follows a chi-squared distribution with $(K - 1)p(p + 1)/2$ degrees of freedom. However for the two sample case, LRT fails when $p > \min(n, m)$ because at least one of $\mathcal{S}_1$ or $\mathcal{S}_2$ will become singular. Bai \etal \cite{Bai2009} and Jiang \etal \citep{Jiang2012} provided asymptotic corrections to the LRT when $n, p \rightarrow \infty$ with $c_n = p/n \rightarrow c \in (0, \infty)$ and proposed
\begin{equation*}
\mathcal{L}^* = \frac{ \mathcal{L} - p \left[ 1 - \left(1 - \frac{n}{p} \right) \log \left( 1 - \frac{p}{n} \right) \right] - \frac{1}{2} \log\left(1 - \frac{p}{n}\right)}{\sqrt{ -2 \left[ \log \left( 1 - \frac{p}{n} \right) - \frac{p}{n} \right] }}
\end{equation*}
which is asymptotically normally distributed under the null hypothesis. 

Another approach for testing equality of covariance matrices is to construct a functional $\mathcal{F}(\Sigma_1, \Sigma_2)$ which will be equal to zero when $\Sigma_1 = \Sigma_2$. Schott \citep{Schott2007} used the squared Frobenius norm of the difference $\Sigma_1 - \Sigma_2$ as the functional to base the test statistic. This method is readily extended to comparing $K$ covariance matrices, with the test statistic $\mathcal{T} = \mathcal{F}_n/\sqrt{\widehat{{\rm var}(\mathcal{F}_n)}}$, where
\begin{gather}
\mathcal{F}_n = \sum \limits_{i < j} \tr \left\{ (\mathcal{S}_i - \mathcal{S}_j)^2 \right\} - (K - 1) \sum \limits_{i = 1}^K \frac{1}{n_i \eta_i} \left[ n_{i} \left( n_{i} - 2 \right) \tr \left( \mathcal{S}_{i}^{2} \right) + n_{i}^{2} \left\{ \tr \left( \mathcal{S}_{i} \right) \right\}^{2} \right],
\end{gather}
and $\eta_i = (n_i + 2)(n_i - 1)$. When $p/n_i \rightarrow c_i \in [0, \infty)$, $\mathcal{T}$ is asymptotically normal under the null hypothesis. Srivastava \etal \citep{SrivastavaCov2010, SrivastavaCov2013, SrivastavaCov2014} developed test statistics using similar rationale but replacing normality assumption with constraints on moments of first four orders. Relaxing the direct relationship between $p$ and $n$, Li and Chen \citep{Li2012} proposed a test statistic by using the $\tr \{(\Sigma_1 - \Sigma_2)^2\}$ as the functional. The test statistic was constructed using U-statistics of the form $\{n(n-1)\}^{-1} \sum_{i < j} (\xvec_i^{\top} \xvec_j)^2$ to estimate $\tr\{\Sigma_1^2\}$ and so on. The leave-out cross-products in the proposed test statistic is similar in spirit to the variance estimate in $T_{CQ}$ \eqref{eqn:chenqin}. Assumptions for the test statistic are similar to \ref{itm:cq3} and \ref{itm:cq4}. 

Covariance matrix estimation is an exciting field which direct applications in graphical network models. Most theory of regularization based sparse precision matrices is based on Gaussian distributions. Extending such estimation to distributions such as Dirichlet-Multinomial or multivariate Poisson where the covariance matrix is parameterized through the mean is very challenging. Hypothesis tests for covariance matrices have primarily been developed by studying the asymptotic properties of traditional test statistics. As seen in Section \ref{sec:meanvt}, random projection methods show good promise in mean vector testing. Using random projections for covariance matrices is an interesting question that is an active area of research. If $R \in \mathbb{R}^{k \times p}$ is an orthogonal random matrix, then projecting the data using $R$ preserves the hypotheses of sphericity and identity in equation \eqref{eqn:covmathyps}. The hypotheses conditional on the random projections will be
\begin{align}
\mbox{Sphericity: } \hspace{1cm} & H_0 : R \Sigma R^{\top} = \sigma^2 R \mathcal{I} R^{\top} (= \sigma^2 \mathcal{I}) \hspace{5mm} \mbox{ vs. } \hspace{5mm} H_A : \Sigma \neq \sigma^2 R \mathcal{I} R^{\top} (= \sigma^2 \mathcal{I}) \nonumber \\
\mbox{Identity: } \hspace{1cm} & H_0 : \Sigma = R \mathcal{I} R^{\top} (= \mathcal{I}) \hspace{5mm} \mbox{ vs. } \hspace{5mm} H_A : \Sigma \neq R \mathcal{I} R^{\top} (= \mathcal{I}).
\label{eqn:covmathyps3}
\end{align}
Theoretical properties of such tests are an active area of research.

\section{Discrete multivariate models}
\label{sec:discretemodels}
Multivariate count data occur frequently in genomics and text mining. In high-throughput genomic experiments such as RNA-Seq (Wang \etal \citep{Zhong2009}), data is reported as the number of reads aligned to the genes in a reference genome. In text mining (Blei \etal \citep{Blei2003}), the number of occurrences of a {\it dictionary} of words in a {\it library} of books is counted to study patterns of keywords and topics. In metagenomics (Holmes \etal \citep{Holmes2012}), abundances of bacterial species in samples is studied by recording the counts of reads assigned to different bacterial species. In all data sets, the data matrix consists of non-negative integer counts. Analyzing multivariate discrete data can be addressed two ways. The absolute counts can be modeled using discrete probability models or the data can be transformed (e.g. using relative abundances instead of absolute counts) and use continuous probability models such as Gaussian, etc. The research community is still divided in opinion on the loss of information due to this transformation (McMurdie and Holmes \citep{McMurdie2014}) or the lack thereof. Transforming the variables will enable us to use hypothesis testing tools presented in Section \ref{sec:meanvt}. In this section, we will look at some discrete multivariate models.

\subsection{Multinomial distribution}
The Multinomial distribution is the most commonly used multivariate discrete model, extending the univariate binomial distribution to multiple dimensions. For $p \geq 2$, the multinomial distribution is parameterized by a probability vector $\pivec = (\pi_1, \ldots, \pi_p)$ with $\pi_1 + \ldots + \pi_p = 1$ and the total count $N \in \mathbb{Z}_+$. The probability mass function of $\xvec \sim {\rm Mult} \left(N, \pivec \right)$ is given by
\begin{equation}
P(\xvec = \mathbf{x}) = P \left(X_1 = x_1, \ldots, X_p = x_p \right) = \frac{ N!}{x_1! \cdots x_p!} \pi_1^{x_1} \cdots \pi_p^{x_p},
\end{equation}
for all $\mathbf{x} \in \mathbb{Z}_+^p$ such that $x_1 + \ldots + x_p = N$. An alternative representation of the multinomial distribution can be obtained using independent Poisson random variables. Consider $p$ independent Poisson random variables, $X_k \sim {\rm Pois}(\lambda_k), k = 1, \ldots, p$. Then the vector $(X_1, \ldots, X_p)$, conditional on $\msum_{k = 1}^p X_k = N$, follows a multinomial distribution with probability parameter $\pivec = (\lambda_1, \ldots, \lambda_p)/(\lambda_1 + \ldots + \lambda_p)$. The re-parameterization using Poisson variable is scale invariant, i.e. the same multinomial distribution is obtained when $X_k \sim {\rm Pois} \left( s \lambda_k \right)$ for all $s > 0$. Levin \citep{Levin1981} provide a very simple expression for the cumulative distribution function using this property,
\begin{equation}
F_{\xvec}(a_1, \ldots, a_p) = P \left( X_1 \leq a_1, \ldots, X_p \leq a_p \right) = \frac{N!}{s^N e^{-s}} \left\{ \prod_{k = 1}^p P \left(Y_k \leq a_k \right) \right\} P(S = N),
\label{eqn:multinomcdf}
\end{equation}
where $s > 0$ is any positive number, $X_k \sim {\rm Pois} \left( s \pi_k \right)$ and $S = Y_1^* + \ldots + Y_p^*$ where $Y_k^*$ is a truncated Poisson variable, $Y_k \sim {\rm Pois}(s \pi_k; \{0, \ldots, a_k \})$. This alternative formulation and equation \eqref{eqn:multinomcdf} reduce the computational cost of calculating the CDF significantly. Using the mass function, the calculation would include doing a comprehensive search in the sample space $\{\mathbf{X}: X_1 + \ldots + X_p = N \}$, which has a computational cost of exponential order with respect to $p$.

The first two moments are functions of $\pivec$, given by $\mathbb{E}(\xvec) = N \pivec$ and ${\rm var}(\xvec) = N \left\{ {\rm diag}\left( \pivec + \pivec^2 \right) - \pivec \pivec^{\top} \right\}$. The constraint on the total sum implies the variables are always negatively correlated, with $\cov \left(X_i, X_j \right) = -N \pi_i \pi_j$. Parameter estimation for multinomial distributions is a well studied. Using the added constraint $\pi_1 + \ldots + \pi_p = 1$, the maximum likelihood estimates can be easily derived as
\begin{equation}
\widehat{\pi}_k = \frac{X_k}{N}, \hspace{1cm} k = 1, \ldots, p.
\label{eqn:multinommle}
\end{equation}
Starting with the works by Rao \citep{Rao2016, Rao1957} wherein consistency and asymptotic properties of the maximum likelihood estimator have been established, several extensions have been developed. When $\pivec$ is restricted to a convex region in the parameter space, Barmi and Dykstra \citep{Barmi1994} developed an iterative estimation method based on a primal-dual formulation of the problem. Jewell and Kalbfleisch \citep{Jewell2004} developed estimators when the multinomial parameters are ordered, i.e. $\pi_1 \leq \pi_2 \leq \ldots \leq \pi_p$. Leonard \citep{Leonard1977} provided a Bayesian approach to parameter estimation by imposing a Dirichlet prior on the probability vector and derived the Bayesian estimates under a quadratic loss function. 
 
When comparing two multinomial populations, $\xvec \sim {\rm Mult}(\pivec_{X})$ and $\yvec \sim {\rm Mult}(\pivec_{Y})$, the hypothesis of interest is
\begin{equation}
H_0 : \pivec_X = \pivec_Y \hspace{1cm} \mbox{ vs. } \hspace{1cm} H_A : \pivec_X \neq \pivec_Y.
\label{eqn:multhyptest}
\end{equation}
Unlike the hypothesis tests in Section \ref{sec:meanvt}, we do not require replicates of the count vectors to construct the test statistic and study its asymptotic properties. Instead, sample sizes for \ref{eqn:multhyptest} are $n = \msum_{k = 1}^p X_k$ and $m = \msum_{k = 1}^p Y_k$.   Traditional tests include the Pearson chi-squared test and the likelihood ratio test, 
\begin{equation}
T_{Pearson} = \msum_{k = 1}^p \frac{ \left(X_k - \widehat{X}_k \right)^2}{\widehat{X}_k} + \frac{ \left(Y_k - \widehat{Y}_k \right)^2}{\widehat{Y}_k}, 
T_{LRT} = \msum_{k = 1}^p \left\{ X_k \log \left(\frac{\widehat{\pi}_k}{\widehat{\pi}_{Xk}} \right) + Y_k \log \left(\frac{\widehat{\pi}_k}{\widehat{\pi}_{Yk}} \right)   \right\},
\end{equation}
where $\widehat{\pi}_k = (X_k + Y_k)/(n + m), \widehat{\pi}_{Xk} = X_k/n, \widehat{\pi}_{Yk} = Y_k/m, \widehat{X}_k = n \widehat{\pi}_k$ and $\widehat{Y}_k = m \widehat{\pi}_k$. Asymptotically, the tests follow a chi-squared distribution with $p$ degrees of freedom under $H_0$. When $p$ is fixed, Hoeffding \citep{Hoeffding1965} provided asymptotically optimal tests for \eqref{eqn:multhyptest}. Furthermore, he also provided conditions under which $T_{LRT}$ has superior performance compared to $T_{Pearson}$. Morris \citep{Morris1975} provided a general framework for deriving the limiting distributions of any general sums of the form 
\begin{equation*}
\mathcal{S}_p = \mathop{\sum}_{k = 1}^p f_k(X_k)
\end{equation*}
when $\{f_k, k = 1, \ldots, p \}$ are polynomials of bounded degree, which generalize $T_{Pearson}$ and $T_{LRT}$. For a comprehensive review of tests, refer to \citep{Balakrishnan2017} and the references therein.

Distributional properties of these tests hold valid when all the counts are large, i.e. $X_k > 0$ and $Y_k > 0$ and number of categories $p$ is smaller than $n + m$. When $p$ is larger than $n$ we encounter sparsity. This is because the minimum number of zero elements will be $p - (n + m)$. Results derived by Morris hold when $p$ and $n + m$ both increase. When the data is large and sparse, i.e. $p > n + m$, Zelterman \citep{Zelterman2019} derived the mean and standard deviation of $T_{Pearson}$ and normalized the test statistic to construct an asymptotically normal test statistic. Using the $\ell_1$ norm of difference, $\| \pivec_{X} - \pivec_{Y} \|_1 = \msum_{k = 1}^p |\pi_{Xk} - \pi_{Yk}|$, and the Euclidean norm $\| \pivec_X - \pivec_Y \|_2^2 = \msum_{k = 1}^p (\pi_{Xk} - \pi_{Yk})^2$ Chan \etal \citep{Chan2014} the following functionals to use as test statistics:
\begin{equation}
\mathcal{T}_1 = \msum_{k = 1}^p \frac{ \left(X_k - Y_k \right)^2 - X_k - Y_k}{X_k + Y_k}, \hspace{1cm} 
\mathcal{T}_2 = \msum_{k = 1}^p \left(X_k - Y_k\right)^2 - X_k  - Y_k
\end{equation}
However, the sampling distributions of $\mathcal{T}_1$ and $\mathcal{T}_2$ were not provided. Instead, permutation based cut-off need to be calculated to do inference.  

Studying the asymptotic properties of such functionals, Plunkett and Park \citep{Plunkett2018} constructed a test statistic, given by
\begin{equation}
T_{PP} = \frac{ \msum_{k = 1}^p \left\{ \left( \frac{X_k}{n} - \frac{Y_k}{m} \right)^2 - \frac{X_k}{n} - \frac{Y_k}{m} \right\} }{\sqrt{ \msum_{k = 1}^p \frac{2}{n^2} \left( \widehat{\pi}_{Xk}^2 - \frac{\widehat{\pi}_{Xk}}{n} \right) + \frac{2}{m^2} \left( \widehat{\pi}_{Yk}^2 - \frac{\widehat{\pi}_{Yk}}{m} \right) + \frac{4}{nm} \widehat{\pi}_{Xk} \widehat{\pi}_{Yk}   }}.
\label{eqn:plunkettpark}
\end{equation}
The test statistic was shown to be asymptotically normal under the following conditions:
\begin{enumerate}[label=(PP \Roman*)]
	\item \label{itm:pp1} $\min(n, m) \rightarrow \infty$ and $n/(n + m) \rightarrow c \in (0,1)$.  This condition is the same as \ref{itm:bs2}, \ref{itm:sd2} and \ref{itm:pa2}.
	\item \label{itm:pp2} The probabilities are not {\it concentrated}, i.e.
	\begin{equation*}
	\frac{1}{\|\pivec_X\|_2^2}	\max \limits_k \pi_{Xk}^2 \rightarrow 0 \mbox{  and  } \frac{1}{\|\pivec_Y\|_2^2} \max \limits_k \pi_{Yk}^2 \rightarrow 0 \mbox{ as } p \rightarrow \infty.
	\end{equation*}
	This condition ensures that the number of components with non-zero probabilities is not bounded. For example, we cannot have  $\pivec_X = (1/m, \ldots, 1/m, 0, \ldots, 0)$ where the number of non-zero elements is equal to $m$ because $\max_k \pi_{Xk}^2 = 1/m^2$ and $\| \pivec_X \|_2^2 = 1/m$ resulting in the ratio being equal to 1/m. 
	\item The sample sizes $n$ and $m$ and dimension $p$ are restricted as
	\begin{equation*}
	(n + m) \| \pivec_X + \pivec_Y \|_2^2 \geq \epsilon > 0 \mbox{ for some } \epsilon > 0.
	\end{equation*}
	To better understand this condition, consider $\pivec_X + \pivec_Y = (1/p, \ldots, 1/p)$. Then $(n + m) \| \pivec_X + \pivec_Y \|_2^2 = (n + m)/p$ which implies $p$ can increase at most linearly with respect to $n$. 
	\item Asymptotic normality is valid in the local alternative
	\begin{equation*}
	n^2 \| \pivec_X - \pivec_Y \|_2^2 = O \left( \| \pivec_X + \pivec_Y \|_2^2 \right).
	\end{equation*}
\end{enumerate}

\subsection{Compound Multinomial models}
\label{sec:compmult}

Consider $n$ multivariate count vectors of dimension $p$, $\xvec_1, \ldots, \xvec_n$. Such data commonly arises when multiple samples are collected, e.g. gene expression counts of $p$ genes collected from $n$ specimens. One common criticism of the standard multinomial distribution is that it does not address {\it over-dispersion} in the data. If we consider that the count vectors are \iid from ${\rm Mult}(\pivec)$, we are inadvertently assuming that the population is homogeneous. To account for heterogeneity in the population, it is advised to assume a model with sample-specific parameter,
\begin{equation*}
\xvec_i | \pivec_i \sim {\rm Mult}(\pivec_i), \hspace{1cm} i = 1, \ldots, n.
\end{equation*}
The heterogeneity can further be modeled using a distribution on the $p$-dimensional simplex $\mathcal{S}_p = \{ \pivec \in \mathbb{R} : \pi_1 + \cdots + \pi_p = 1 \}$. In the univariate case, the beta distribution is the natural choice for the distribution on $\mathcal{S}_2$. Extending to $p$ dimensions, the natural extension is the multivariate beta distribution or the Dirichlet distribution. 

The Dirichlet distribution is characterized by a single parameter $\thetavec = (\theta_1, \ldots, \theta_p)$, with density function
\begin{equation*}
f(\pivec; \thetavec) = \frac{ \Gamma(\theta_0)}{\prod \limits_{k = 1}^p \Gamma(\theta_k)} \pi_1^{\theta_1 - 1} \cdots \pi_p^{\theta_p}, \hspace{5mm} \pi_1 + \cdots + \pi_p = 1.
\end{equation*}
where $\theta_0 = \theta_1 + \cdots + \theta_p$ and $\Gamma(\cdot)$ is the gamma function. The compound Dirichlet-Multinomial(DirMult) distribution, constructed by the marginal of $\xvec_i | \pivec_i \sim {\rm Mult}(\pivec_i)$ and $\pivec \sim {\rm Dir}(\thetavec)$ has the density function given by
\begin{equation}
f(\xvec; \thetavec) =  \frac{\Gamma(X_0 + 1) \Gamma(\theta_0)}{\Gamma(X_0 + \theta_0)} \prod \limits_{k = 1}^p \frac{ \Gamma(X_k + \theta_k)}{ \Gamma(X_k + 1) \Gamma(\theta_k)},
\label{eqn:dirmultpdf}
\end{equation}
where $X_0 = X_1 + \cdots + X_p$. The DirMult model was first introduced by Mosimann, who derived the properties of the distribution. The mean and variance of the DirMult distribution are $\mathbb{E}(\xvec) = X_0 \theta_0^{-1} \thetavec$ and ${\rm var}(\xvec) = n \left\{ \theta_0^{-1} {\rm diag}(\thetavec) - \theta_0^{-2} (X_0 + \theta_0)/(1 + \theta_0) \thetavec \thetavec^{\top}  \right\}$. The variance matrix is the sum of a full-rank matrix (diagonal part) and a rank-one matrix. Using the result from Miller \citep{Miller1981}, the precision matrix can be calculated in closed form as
\begin{equation*}
{\rm var}(\xvec)^{-1} = n^{-1} \left\{ \theta_0 {\rm diag}(\thetavec^{-1}) + \frac{ \frac{X_0 + \theta_0}{1 + \theta_0}}{\theta_0^2 - \frac{X_0 + \theta_0}{1 + \theta_0}}  \mathbf{1} \mathbf{1}^{\top}     \right\}.
\end{equation*}

For parameter estimation, the likelihood function of \eqref{eqn:dirmultpdf} does not admit a maximum for $\thetavec$ in closed form. An approximate solution can be obtained using iterative methods such as the Newton-Raphson algorithm. One convenient feature for computation is that the second-order derivative of the log-likelihood function has a closed-form expression for the inverse (Sklar \citep{Sklar2014}). Thus the Newton-Raphson step has a linear computation cost. When $p$ is larger than $X_0$, Danaher \citep{Danaher1988} derived parameter estimates the beta-binomial marginals and established their consistency.

While the density function is known to be {\it globally convex}, maximization can still lead us to a local maxima. A proper initial value specification is essential to have good performance of the estimator. Choice of optimal initial values has been an area of considerable interest, even for the Dirichlet distribution. The challenge lies in the fact that the method of moments (MM) estimator is {\it not unique}. This is because of the scaling in $\mathbb{X}_k = n \theta_k/\theta_0$, which gives both $\widehat{\thetavec}$ and $c \widehat{\thetavec}$ as MM estimates for any $c > 0$. Ronning \citep{Ronning1989} proposed using the same initial value for all elements, $\widehat{\theta}_k = \min \limits_{ij} X_{ij}$. This proposal was based on an observation that the method of moments estimates can lead to Newton-Raphson updates becoming inadmissible, i.e. $\widehat{\theta}_k < 0$ for some $k$. Hariharan \citep{Hariharan1993} have done a comprehensive comparison of the different initial values under several models. However they concluded that none of the methods is uniformly consistent across all the models.

Dirichlet-Multinomial has been applied to study multivariate count data in several applications in biomedical research. In metagenomics, the study of bacterial composition of environmental (biological or ecological) samples, we are interested in modeling the abundance of different species of bacteria in samples. The Dirichlet-Multinomial model is apt for such data because (i) abundances of bacteria are constrained by the total number of bacteria sampled in the specimen and (ii) over-dispersion due to environmental variability is accounted for. Holmes \etal \citep{Holmes2012} used a Dirichlet multinomial mixture model to cluster samples by abundance profile, i.e. the DirMult parameter. Chen and Li \citep{Chen2013} developed a $\ell_1$-penalized parameter estimation for variable selection in the DirMult model. Sun \etal \citep{Sun2018} used the DirMult model to construct a clustering algorithm for single-cell RNA-seq data.

The most celebrated application of DirMult distribution is latent dirichlet allocation (LDA), introduced by Blei \etal \citep{Blei2003}. Developed for text mining for classifying documents by keywords, the model is a hierarchical Bayesian model with three levels. Firstly, the $p$ elements of $\xvec$ represent the {\it words} in the vocabulary. A word is represented as $\xvec = (x_1, \ldots, x_p)$ where $x_k \in \{0, 1\}$ for all $k = 1, \ldots, p$ and $\msum_{k = 1}^p x_k = 1$. A collection of $q$ words represents a {\it topic}, which can be used to classify documents, which will also be a multinomial variable $\mathbf{T} = (t_1, \ldots, t_q)$ with $t_k \in \{0, 1\}$ for all $k = 1, \ldots, q$. The number of topics, $K$, is assumed to be fixed. It should be noted that while the words in the vocabulary are defined and observed, the topic corresponding to a word is a {\it latent} variable. Second, each {\it document} is defined as a sequence of $N$ words, $\mathcal{X} = \{ \xvec_1, \ldots, \xvec_N \}$. The number of words in a document is assumed to have a Poisson distribution ($N \sim {\rm Pois}(\lambda)$) and the topics follow a multinomial distribution with document-specific parameter. And finally, a {\it corpus} is defined as a collection of $M$ documents, $ \mathcal{D}_N = \{\mathcal{X}_1, \ldots, \mathcal{X}_M \}$. 

The LDA model is parameterized as follows. Each corpus is characterized by the probability of its keywords $\thetavec_m$, $\mathbf{T} \sim {\rm Mult}(\thetavec_m)$. The probability parameters are assumed to be following a Dirichlet distribution, $\thetavec_m \sim {\rm Dir}(\boldsymbol{\alpha}), m = 1, \ldots, M$. Conditional on the latent topics $\mathbf{T}$, $\pi_{kt} = P \left( X_k = 1 | T_t = 1 \right)$ denotes the probability that $k^{\rm th}$ word in the vocabulary is observed, provided the word describes the topic. The collection of all such probabilities is parameterized as a $p \times q$ matrix $\boldsymbol{\Pi} = (\pi_{kt} : k = 1, \ldots, p; t = 1, \ldots, q)$. Using these components, the complete likelihood can be written as
\begin{align}
P \left( \mathbf{D} | \boldsymbol{\alpha}, \boldsymbol{\Pi} \right) & = \prod \limits_{d = 1}^D \quad \mathop{\int} \limits_{\mathcal{S}_p} f(\thetavec_d | \boldsymbol{\alpha}) P \left( \mathcal{X}_d | \thetavec, \boldsymbol{\Pi} \right) \, d \thetavec_d \nonumber \\
& = \prod \limits_{d = 1}^D \quad \mathop{\int} \limits_{\mathcal{S}_p} f(\thetavec_d | \boldsymbol{\alpha}) \left\{ \prod \limits_{n = 1}^{N_d}  g \left( \mathbf{T}_n | \thetavec_d \right)  h \left( \xvec_{n} | \mathbf{T}_n, \boldsymbol{\Pi} \right) \right\} d \thetavec_d.
\label{eqn:ldamodel}
\end{align}
In this model, $f(\cdot)$ is the Dirichlet density function, $g(\cdot)$ is the multinomial mass function and $h(\cdot)$ is obtained from $\boldsymbol{\Pi}$. Parameter estimation is done by maximizing the likelihood using expectation-maximization (EM) algorithm by conditioning on the latent keywords. 

Major focus on LDA research has been on developing faster algorithms (Hoffman \etal \citep{Hoffman2010}) to be able to analyze larger corpora with large number of documents. Mimno \etal \citep{Mimno2012} considered sparsity in the model from the Gibbs sampling perspective to improve the efficiency of the algorithm. However most of the research has been from a {\it machine learning} and estimation perspective. Statistical properties of the estimators, which could be of potential interest for developing hypothesis tests, have not been established. One potential problem of interest could be comparing the Dirichlet parameters of two corpora,
\begin{equation}
H_0 : \boldsymbol{\alpha}_1 = \boldsymbol{\alpha}_2 \hspace{1cm} \mbox{ vs. } \hspace{1cm} H_A : \boldsymbol{\alpha}_1 \neq \boldsymbol{\alpha}_2.
\label{eqn:ldamodel2}
\end{equation} 
In computer science literature, the focus has been on developing methods for efficient analysis of corpora with large number of documents. Sample size is known to affect accuracy of the allocation (Crossley \etal \citep{Crossley2017}). A large $p$ small $n$ problem in this context would be efficient classification of small number of documents (small N) with a large vocabulary (large p). Understanding the efficiency of LDA in such large $p$ small $n$ scenarios is an open area of research.

\subsection{Other distributions}

The Dirichlet-Multinomial is a natural extension to the univariate beta-binomial distribution, which are the marginals of the DirMult distribution. This observation arises the following question: can we develop multivariate count distributions with known marginals? The theoretical answer to this question is to use Sklar's theorem (Nelson \citep{NelsonBook}) and construct a copula to model the joint distribution. However parametric inference such as hypothesis testing is very tedious and sometimes intractable when using copula models. In this section, we shall look at some multivariate extensions to known univariate distributions which have useful parameterizations and are easy to do inference. 

\subsubsection{Bernoulli distribution}
One of the earliest generalizations of the Bernoulli distribution using a parametric approach was developed by Teugels \citep{Teugels1990}. Using the moments of all orders $k = 1, \ldots, p$, the moment generating function of multivariate Bernoulli was constructed. They also provided an extension to the multivariate binomial distribution using the sum of independent Bernoulli variables. Using the joint probabilities, Dai \etal \citep{Dai2013} proposed a multivariate Bernoulli distribution which has an analytical form of the mass function. Before generalizing the multivariate Bernoulli distribution, consider the case where elements of the variable $\xvec = (X_1, \ldots, X_p)$ are independent with $X_k \sim {\rm Ber}(\pi_k), k = 1, \ldots, p$. Then the joint probability of $\xvec = \mathbf{x}$ is given by
\begin{equation*}
P \left(\xvec = \mathbf{x} \right) = \prod \limits_{k = 1}^p P(X_k = x_k) = \prod \limits_{k = 1}^p \pi_k^{x_k} \left(1 - \pi_k \right)^{1 - x_k}.
\end{equation*}
When the variables are dependent, the joint probability cannot be factored into the product of marginals. Using the joint probabilities, the mass function can defined as
\begin{equation}
P \left(X_1 = x_1, \ldots, X_p \right) = \pi_{00\ldots0}^{\prod\limits_{k = 1}^p (1 - x_k)} \times \pi_{10\ldots0}^{x_1 \prod\limits_{k = 2}^p (1 - x_k)} \times \cdots \times \pi_{11\ldots1}^{\prod\limits_{k = 1}^p x_k}, 
\label{eqn:multbernoulli}
\end{equation}
where $\pi_{00\ldots0} = P(X_1 = 0, \ldots, X_p = 0)$ and so on. The marginals of $\xvec$ are Bernoulli with cumulative probability,
\begin{equation*}
X_k \sim {\rm Ber}(\pi_k), \hspace{1cm} \pi_k = \sum \limits_{i \neq k: a_i = 0, 1} \pi_{a_1\ldots a_{k-1}1a_{k+1}\ldots a_p}.
\end{equation*}
Using this formulation, they computed the moments and also calculate maximum likelihood estimates using Newton-Raphson algorithm. However the main drawback is the dimension of the parameter space. To define the multivariate Bernoulli mass function, we require a total of $2^p - 1$ parameters, which can be computationally infeasible for higher dimensions. 

\subsubsection{Binomial distribution}
The bivariate binomial distribution (BBD) was first introduced by Aitken and Gonin \citep{Aitken1936} in the context of analysis $2 \times 2$ contingency tables when the two outcomes are not independent. Several extensions have been provided since, including work by Krishnamoorthy \citep{Krishnamoorthy1951} who derived the properties of BBD by extending the moment-generating function from the independent case to dependent variables. Hudson and Tucker \citep{Hudson1986} established limit theorems for BBD expressing them as sums of independent multivariate Bernoulli variables. Several other researchers have discussed the properties of BBD. For a recent list of all publications, please refer to Biswas and Hwang \citep{Biswas2002} and the references therein. The multivariate binomial distribution (MBD) also suffers from the same {\it curse of dimensionality} as the Bernoulli distribution. The total number of parameters required to define the $p$-dimensional distribution is equal to $2^p - 1$.

The multivariate binomial distribution poses several questions that still need to be answered. For instance, it would of interest to simplify the distribution for a restricted parameter set. For instance, if we assume only $k$-fold interactions are feasible, then the model can be reduced to have $2^k - 1$ parameters. The generalized additive and multiplicative binomial distribution models proposed by Altham \citep{Altham1978} can serve as motivation for building such reduced models. MBD can also be used to model several data sets in genomics. For instance when studying epigenomic modifications such as DNA methylation, co-methylation (mutual methylation of pairs of genes) is actively studied for understanding their association with different phenotypes (outcomes). MBD can be used to model the joint probability of methylation of pairs of genes. However the major bottleneck that needs to be solved first is the computational complexity. Currently, there are no existing tools to compute and model MBD. With improved computational capabilities, this task should be accomplished easily. 


\subsubsection{Poisson distribution}

Constructing a multivariate Poisson distribution whose marginals are univariate Poisson variables is fairly easy. Consider the bivariate case. If $Z_k \sim {\rm Pois}(\lambda_k), k = 1, 2, 3$ are independent Poisson variables, then $\xvec = (X_1, X_2)$ defined as
\begin{equation*}
X_1 = Z_1 + Z_3, \hspace{1cm} X_2 = Z_2 + Z_3
\end{equation*}
gives a bivariate distribution with Poisson marginals, $X_1 \sim {\rm Pois}(\lambda_1 + \lambda_3)$ and $X_2 \sim {\rm Pois}(\lambda_2 + \lambda_3)$. The joint mass function can be expressed as
\begin{align}
P(X_1 = x_1, X_2 = x_2) & = \sum \limits_{z = 0}^{\min(x_1, x_2)} P(Z_1 = x_1 - z, Z_2 = x_2 - z, Z_3 = z) \nonumber \\
& = e^{-(\lambda_1 + \lambda_2 + \lambda_3)} \sum \limits_{z = 0}^{\min(x_1, x_2)}  \frac{ \lambda_1^{x_1 - z}}{(x_1 - z)!} \frac{\lambda_2^{x_2 - z}}{(x_2 - z)!} \frac{\lambda_3^z}{z!}.
\label{eqn:multpoisson}
\end{align}
Extending to $p$ dimensions, the multivariate Poisson is defined through the latent $Z_k$'s as
\begin{equation}
X_k = Z_{kk} + \sum \limits_{j \neq k} Z_{kj}, \hspace{1cm} k = 1, \ldots, p,
\label{eqn:multpois2}
\end{equation}
where $Z_{jk} \sim {\rm Pois}(\lambda_{jk})$. Expressing the latent variables in matrix form $(Z_{jk})_{j,k = 1, \ldots, p}$, defining $\xvec$ requires $p(p + 1)/2$ independent latent components. The mass function can be expressed as $p(p - 1)/2$ summations and is computationally intensive for even moderate values of $p$. A more general form of the multivariate Poisson requires $2^p - 1$ latent components and is infeasible to express as in equation \eqref{eqn:multpois2}. The following trivariate Poisson should serve as a basic overview of the idea:
\begin{align}
X_1 & = Z_{1} + Z_{12} + Z_{13} + Z_{123}, \nonumber \\
X_2 & = Z_{2} + Z_{12} + Z_{23} + Z_{123}, \nonumber \\
X_3 & = Z_{3} + Z_{13} + Z_{23} + Z_{123} 
\label{eqn:multpoisexample}
\end{align}

The main drawback with this formulation of multivariate Poisson distribution is its restrictive dependence structure. In the bivariate case, the correlation between $X_1$ and $X_2$ is given by
\begin{equation*}
{\rm cor} \left( X_1, X_2 \right) = \frac{ \lambda_{3}}{\sqrt{\lambda_1 + \lambda_3} \sqrt{ \lambda_2 + \lambda_3} },
\end{equation*}
which is always positive. Extending the distribution to a larger class of correlation structures, Shin and Pasupathy \citep{Shin2010} proposed using the {\bf nor}mal {\bf t}o {\bf a}nything (NORTA) algorithm \citep{Cario1997} for random number generation from multivariate Poisson with negative correlations. They define the iterative procedure for generating bivariate Poisson variables with correlation $\rho$ as follows. Let $U_1, U_2, U_3 \sim U(0, 1)$ be \iid variables. A bivariate Poisson distribution with marginals $X_1 \sim {\rm Pois}(\lambda_1)$ and $X_2 \sim {\rm Pois}(\lambda_2)$ can be obtained using
\begin{equation}
X_1 = F^{-1}_{\lambda_1 - \lambda*}(U_1) + F^{-1}_{\lambda^*}(U_3), \hspace{5mm} 
X_2 = \left\{  \begin{array}{ll}F^{-1}_{\lambda_2 - \lambda_2 \lambda^*/\lambda_1}(U_2) + F^{-1}_{\lambda_2 \lambda^*/\lambda_1}(U_3) & \hspace{5mm} \mbox{ if } \rho > 0  \\
\\
F^{-1}_{\lambda_2 - \lambda_2 \lambda^*/\lambda_1}(U_2) + F^{-1}_{\lambda_2 \lambda^*/\lambda_1}(1 - U_3) & \hspace{5mm} \mbox{ if } \rho < 0 \end{array} \right.,
\label{eqn:multpoistrex}
\end{equation}
where $F^{-1}_{\lambda}(x) = \inf \{y : F_{\lambda}(x) \geq y\}$ is the inverse Poisson cumulative distribution function with parameter $\lambda$. The parameter $\lambda^* \in [0, \lambda_1]$ assuming $\lambda_1 \leq \lambda_2$. If $\lambda_1 \geq \lambda_2$, $X_1$ and $X_2$ can be inter-changed. While this formulation gives a method for generating random samples from bivariate Poisson variables with negative correlations, it is unusable for inference as the likelihood function is not available. Obtaining the likelihood function for the bivariate case using \eqref{eqn:multpoistrex} and parameter estimation using the derived likeliho0d are a few open problems in using this construction of multivariate Poisson variables. 

Karlis \citep{Karlis2010} developed another approach to characterize multivariate Poisson random variables by compounding  independent Poisson components through a multivariate distribution on their parameters. If $\boldsymbol{\lambda} = (\lambda_1, \ldots, \lambda_p) \sim \mathcal{G}(\Theta)$ is a multivariate distribution, then dependence structure on $\xvec$ can be imposed by taking the a mixture of independent Poisson distributions with $\mathcal{G}$,
\begin{equation}
P(\xvec = \boldsymbol{x} | \Theta) = \int \limits_{\mathbb{R}^n_{+}} \prod \limits_{k = 1}^p e^{-\lambda_p} \frac{ \lambda_p^{x_p}}{x_p!} g(\boldsymbol{\lambda}; \Theta) \, d \lambda_1 \ldots d \lambda_p
\label{eqn:multpois3}
\end{equation}
A popular choice for $\mathcal{G}$ is the log-normal distribution, since the distribution should be defined on $\mathbb{R}^p_+$. This formulation has two advantages. Firstly, the covariance structure on $\boldsymbol{\lambda}$ will impart a dependence structure on $\xvec$. Secondly, the mixture model ensures that the variablesl of $X_k$ is greater than $\lambda_k$ for all components, thereby addressing issues of over-dispersions. For more details, readers may refer to Inouye \etal \citep{Inouye2018} and the references therein for more papers published studying the multivariate Poisson distribution. 

Multivariate Poisson distributions are fairly new and have a lot of problems that need to be addressed. The framework for hypothesis testing is not extensively developed. There is very limited literature in this regard. For example, Stern \citep{Stern2002} developed a test for the bivariate Poisson model in \ref{eqn:multpoisson} testing for $H_0 : \lambda_3 = 0$ versus $H_A : \lambda_3 \neq 0$ using a Bayesian significance test. Testing hypotheses comparing two or more multivariate Poisson families is not addressed. High dimensional tools for multivariate Poisson are extremely hard to develop due to the exponential computation cost: $2^p - 1$ latent variables required to define the distribution. Restricted models, such as using only pairwise correlations in \eqref{eqn:multpoisexample}, have a quadratic computation cost and are easier to study. These could potentially be a good starting point for studying the complete model. 

%

\section{Conclusion}
\label{sec:conclusion}
High dimensional inference is a very exciting field of statistics with many theoretical challenges and practical uses. Availability of large-scale and high-dimensional data is increasing leaps and bounds.  Conducting large-scale analysis has become practical with the availability of high performance computing facilities. There is an urgent need to develop statistical tools that can tackle these large dimensional data sets efficiently and accurately. Statistical methodology and computational tools need to progress in conjuction with each other, leaving the onus on statisticians to develop more accurate methods for estimation and inference. 

In this chapter, we have addressed three areas of high dimensional inference that are being actively developed. Hypothesis tests for the population mean is one of the more standard inference problems, which has been well studied in high dimensions. We looked at the two main approaches - asymptotics-based tests and random projection based tests have been presented. The asymptotics based tests have been fairly well-studied in comparison to the random projection based tests. Projections into lower-dimensional spaces using random matrices is an active area of research in mean vector testing. We should consider other methods for dimension reduction to study their use in high dimensional inference. Convolutional neural networks (CNN) \citep{Goodfellow2016}, which are commonly used in deep learning, is another exciting dimension reduction technique that is currently not used for high dimensional inference. 

Sparse covariance matrix estimation has found practical use in understanding the graphical network structure of variables in high dimensions. We looked at different approaches to construct the regularization and the computational tools developed for optimization. While Gaussianity of variables is commonly assumed in sparse precision matrix estimation due to its properties, extension to non-Gaussian distributions is to be studied. We have looked at hypothesis testing for comparing two or more covariance matrices in the high dimensional setting. One approach we can identify that is lacking is the use of random projections in covariance matrix testing. This poses an interesting challenge to see the versatility of random projections in high dimensional inference. 

Finally, we looked at development of discrete multivariate models and the challenges therein. Only two distributions have been extensively studied - multinomial and Dirichlet-multinomial.  We looked at high-dimensional hypothesis tests for the multinomial parameters. The hierarchical models and sparse regression models for the Dirichlet-multinomial distribution are also well studied. However a lot of work needs to be done for other distributions. The theoretical developments in multivariate Bernoulli models need to be supplemented with computational tools for estimation and inference. A generalized multivariate Poisson distribution needs to be developed, which can lead to potential extensions such as multivariate Poisson-gamma mixtures.

\bibliographystyle{apalike}
\bibliography{../hos_ayyala_ver5}

\end{document}